\documentclass[]{interact}
\usepackage{xcolor}
\usepackage[ruled,vlined]{algorithm2e}
\usepackage{epstopdf}
\usepackage[caption=false]{subfig}

\usepackage[numbers,sort&compress]{natbib}
\bibpunct[, ]{[}{]}{,}{n}{,}{,}
\makeatletter
\def\NAT@def@citea{\def\@citea{\NAT@separator}}
\makeatother

\theoremstyle{plain}
\newtheorem{theorem}{Theorem}[section]

\theoremstyle{definition}
\newtheorem{definition}[theorem]{Definition}

\theoremstyle{remark}

\begin{document}

\articletype{ARTICLE TEMPLATE}

\title{Optimal control of atmospheric pollution because of urban traffic flow by means of Stackelberg strategies}

\author{
\name{N.~Garc\'\i a-Chan\textsuperscript{a}\thanks{CONTACT N.~Garc\'\i a-Chan. Email: nestor.gchan@academicos.udg.mx},  L.J.~Alvarez-V\'azquez\textsuperscript{b}, A.~Mart\'inez \textsuperscript{b}, M.E.~V\'azquez-M\'endez\textsuperscript{c}}
\affil{\textsuperscript{a}Depto.~F\'\i sica, Universidad de Guadalajara, C.U.~Ciencias Exactas e Ingenier\'\i as,\\ 44430 Guadalajara, Mexico; \textsuperscript{b}Depto.~Matem\'atica Aplicada II, Universidade de Vigo, E.I.~Telecomunicaci\'on, 36310 Vigo, Spain; \textsuperscript{c} Depto.~Matem\'atica Aplicada, Universidade de Santiago de Compostela, E.P.S., 27002 Lugo, Spain}
}

\maketitle

\begin{abstract}
Two major problems in modern cities are air contamination and road congestion. They are closely related and present a similar origin: traffic flow. To face these problems, local governments impose traffic restrictions to prevent the entry of vehicles into sensitive areas, with the final aim of dropping down air pollution levels. However, these restrictions force drivers to look for alternative routes that usually generate congestions, implying both longer travel times and higher levels of air pollution.
In this work, combining optimal control of partial differential equations and computational modelling, we formulate a multi-objective control problem with air pollution and drivers' travel time as objectives and look for its optimal solutions in the sense of Stackelberg. In this problem, local government (the leader) implements traffic restrictions meanwhile the set of drivers (the follower) acts choosing travel preferences against leader constraints. Numerically, the discretized problem is solved by combining genetic-elitist algorithms and interior-point methods, and computational results for a realistic case posed in the Guadalajara Metropolitan Area (Mexico) are shown.
\end{abstract}

\begin{keywords}
Multi-objective,   Numerical simulation, Optimal control, Stackelberg solution, Traffic related air pollution
\end{keywords}

\begin{amscode}  
90C29; 90B50; 49J20 
\end{amscode}

\section{Introduction}

Growth and expansion of major cities have originated, as an undesirable side effect, the critical augmentation of two closely related environmental problems: air pollution and traffic congestions, whose main factor can be considered urban traffic. Regarding the first problem, urban atmospheric contamination is highly subordinate to vehicular emissions (carbon oxides, nitrogen oxides and so on), but concentration levels of such pollutants depend also on other external factors such as, among others, wind or humidity. With respect to the second problem, main negative consequences are related to the increase in the necessary time for its residents to carry out their daily moves, with the resulting discomfort associated, for instance, to excessive fuel consumption, delays and noise pollution.

To confront this problems, common public policies imposed by the local governments are related to traffic restrictions at the intersections of the urban road network. With these restrictions they prevent the entry of vehicles into sensitive areas (normally the city center) with the aim of bringing down the air pollution concentration. However, these limitations force drivers to choose other road preferences to reach their destiny, inducing traffic congestions. Thus, contrary to expectations, these traffic congestions can increase pollutants concentrations and present a negative impact on drivers with longer travel time.

The reduction of air pollution levels by the traffic restrictions and their consequent change in drivers' preferences, is nowadays a controversial topic. Recent studies show that the impact of traffic restrictions (and other traffic management strategies) on air pollution levels is moderately successful in low emissions zones of some European cities. Meanwhile, the lack of data and the complexity of epidemiology studies made harder the detection and identification of traffic-related health impact on inhabitants by exposure to noise, stress and air pollution (see \cite{york} and the references therein).
However, studies also shown that this impact on air pollution can be greater and that society is aware of the need of these traffic restrictions policies. In \cite{invernizzi} the concentration of black carbon was compared with the concentration of particulate matter (PM) in three zones of Milan (Italy): without traffic restrictions, with traffic restrictions and with pedestrians only. This data analysis showed (roughly speaking) that the concentration levels of black carbon drop down from traffic areas to pedestrian areas meanwhile, the PM concentration does it but in a more moderate way. In \cite{pestana} an inquest to the inhabitants of Lisbon (Portugal) showed that they are willing to accept charges for vehicular congestions with the aim of a better quality of life.     

Those studies \cite{york, invernizzi} are mainly empirical, and {\it a priori} evaluation or even the certainty of dropping down the air pollution levels by traffic restrictions are out of their point of view. Therefore, a suitable combination of mathematical models, numerical simulation and optimal control techniques are an important tool for estimating and optimizing the impact of traffic restrictions and drivers' preferences on the air pollution levels and the drivers' travel time. Moreover, these {\it a priori} estimates and minimization results could be employed as a factor to change the viewpoint of city inhabitants, making them agree to these traffic management policies.       

In this context, partial differential equations models are usually employed both in the analysis of urban traffic flow in road networks \cite{holden, coclite, garavello, garavello2, goettlich2015, goatin2016} and in the investigation of atmospheric pollution \cite{alvarez-vazquez2015, garcia-chan1, skyba2, stockie, orun}. Nevertheless, the compounding of both topics has been much less addressed (we can mention, for instance, \cite{Parra-Guevara2003, goettlich2011, Berro, garcia-chan2, Piccoli}), and is usually based on the assumption of a previous knowledge of the vehicular flow, which is not adapted to analyze the management of a road network that may be optimal with respect to travel times and contamination levels. 

The authors have addressed this topic in a series of recent works with a progressive complexity. So, in \cite{alvarez-vazquez2017} a new methodology that couples a 1D model for vehicular flow with a 2D model for pollutant dispersion was proposed, in order to estimate the air pollution related to traffic flow. In \cite{alvarez-vazquez2018} an optimal control problem related to the expansion of an existing urban road network with an environmental perspective was formulated and solved. Finally, in \cite{vazquez-mendez2019} a multi-objective optimal control problem -where the air pollution and the travel time were the objectives, the drivers' preferences were the controls, and the traffic restrictions were fixed- was solved from a cooperative point of view, that is, its Pareto front was obtained using a genetic algorithm.     

Thus, within the framework of the optimal control of partial differential equations, the current work represents a step forward of the authors in the same direction, considering now a non-cooperative, hierarchical point of view, that is, a Stackelberg strategy \cite{stackelberg}. Stackelberg strategies are commonly applied in economy \cite{julien}, and the authors have previously applied them in the optimal management of a wastewater system \cite{alvarez-vazquez2015_2}. Therefore, in the present context of minimizing the urban air pollution levels and the drivers' travel time \cite{vazquez-mendez2019}, the existence of a hierarchical relation between the local government (denoted as the leader) and the set of drivers (denoted as the follower) is assumed. Then a bi-level multi-objective optimal control problem is formulated and its Stackelberg solution is formally defined (Section 3). A complete discretization of the cost functionals and a combination of an interior-point method \cite{waltz} with a genetic-elitist algorithm \cite{goldberg, deb} is proposed to solve this bi-level problem. This combination is carry out using adjoint state techniques \cite{marchuk}, where the pollutant objective functional is written in an alternative, simpler way (that considers the adjoint state and the pollutant emissions both only evaluated on the road network instead of the pollution state evaluated in the whole urban domain), with the aim of reducing the computational cost (Section 4). Finally, some numerical experiences for a real-world case posed in one of the largest metropolitan areas in Mexico (the Guadalajara Metropolitan Area (GMA), with almost five million inhabitants and more than two million vehicles) are presented (Section 5), and several concluding remarks are derived (Section 6).

\section{Mathematical modelling}

\begin{figure}
\begin{center}
\includegraphics[width=8cm]{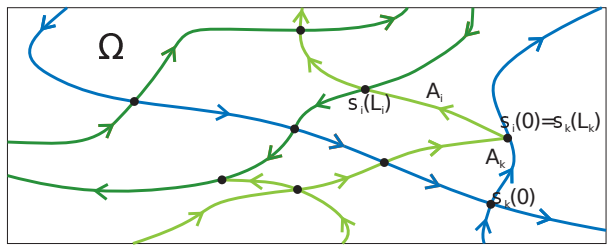}
\end{center}
\caption{Scheme of a typical domain $\Omega$ corresponding to a city with a road network.}
\label{Fig1}
\end{figure}

Let $\Omega\subset\mathbb{R}^2$ be a domain representing a city, including an urban road network formed by $N_R$ unidirectional avenues crossing at $N_J$  intersections, and such that each road endpoint is either an intersection or lies on the boundary of $\Omega$ (see a schematic example in Fig.~\ref{Fig1}). 

Each avenue $A_i\subset\Omega, \ i=1,\dots,N_R,$ is represented by an interval $[0,L_i]$ parametrized by means of the horizontal alignment:
\begin{equation}\label{parametrizacion}
\sigma_i : s \in [0,L_i] \subset\mathbb{R} \longrightarrow \sigma_i(s)=(x_i(s),y_i(s)) \in A_i
\end{equation}
where the arc length parameter $s$ preserves the sense of motion on the road. In the following, we denote by ${\cal I}^{in},\,{\cal I}^{out}\subset\{1,\ldots,N_R\}$ the sets of indices designating incoming and outgoing avenues in the network, respectively, and by ${\cal I}_j^{in},\,{\cal I}_j^{out}\subset\{1,\ldots,N_R\}$ the sets of indices designating incoming and outgoing avenues at the intersection $j \in \{1,\ldots,N_J\}$, respectively.

\subsection{Modelling traffic flow}

\begin{figure}[tb]
\begin{center}
\includegraphics[width=8cm]{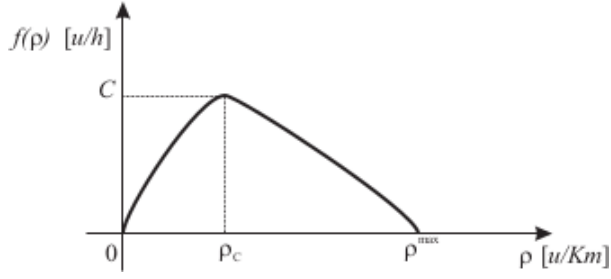}
\end{center}
\caption{Classical static relation showing flow rate $f(\rho)=\rho v$ versus density $\rho$, and depicting maximum density $\rho^{max}$, critical density $\rho_{C}$ and road capacity $C$.}
\label{Fig2}
\end{figure}

In whole road network, traffic flow is modelled by the classical Lighthill-Whitham-Richards (LRW) model coupled with queue terms. Then, we denote by $\rho_i(s,t)\in[0,\rho_i^{max}]$ the density of cars at point $\sigma_i(s)$ of avenue $A_i$ and at time $t\in[0,T]$ (measured in {\it number of cars/km}), where $\rho_i^{max}$ represents the maximum allowed density. The LRW model assumes that the flow rate on each avenue $A_i$ $[\mbox{\it number of cars}/h]$ is given by a function $f_i : [0,\rho_i^{max}]\rightarrow\mathbb{R}$ in terms of the density (i.e., $f_i(\rho_i)=\rho_i v_i$, where $v_i$ $[km/h]$ represents the velocity on the avenue $A_i$). Fundamental diagram $f_i$ is usually known as the {\it static relation} on $A_i$ (see Fig.~\ref{Fig2}), and it must verify the following properties:
\begin{enumerate}
\item $f_i : [0,\rho_i^{max}]\rightarrow\mathbb{R}$ is Lipschitz continuous and concave.
\item $f_i(0)=f_i(\rho_i^{max})=0$.
\item There exists a unique value $\rho_{C_i}\in(0,\rho_i^{max})$ (the so-called {\it critical density}) such that $f_i$ is strictly increasing in $(0,\rho_{C_i})$ and strictly decreasing in $(\rho_{C_i},\rho_i^{max})$ (maximum value $C_i=f_i(\rho_{C_i})$ is known as {\it road capacity}).
\end{enumerate}

Moreover, for all $y\in{\cal I}^{in}$, we define the queue length $q_y(t)\geq 0$ (measured in $\mbox{\it number of cars}$) downstream the avenue $A_y$, where we assume given the desired inflow rate $f^{in}_y(t)$ and the downstream road capacity $C^{in}_y$. We also suppose that, for all $z\in\mathcal{I}^{out}$, the maximum outflow rates $f_z^{out}(t)$ are known.  

So, traffic flow in the whole road network is defined by the solution of the system \cite{garavello, vazquez-mendez2019}: For $i=1,\ldots,N_R$, $y\in{\cal I}^{in}$, $z\in{\cal I}^{out}$, $j=1,\ldots,N_J$, $k\in{\cal I}_j^{in}$, and $l \in{\cal I}_j^{out}$:\\
\begin{subequations}
\begin{align}
\frac{\partial \rho_i}{\partial t} + \frac{\partial f_i(\rho_i)}{\partial s}  =  0  \quad \mbox{ in }(0,L_i)\times(0,T),\label{MT_ED}\\
\rho_i(.,0)=\rho_i^0   \quad \mbox{ in }[0,L_i], \label{MT_IC}\\
f_{k}(\rho_{k}(L_{k},.))=\sum_{l\in{\cal I}_j^{out}}\min\left\{\alpha_{lk}^jD_{k}(\rho_{k}(L_{k},.)),\,\beta_{kl}^j S_{l}(\rho_{l}(0,.))\right\} \quad \mbox{ in }(0,T), \label{MT_CCj_in}\\
f_{l}(\rho_{l}(0,.))=\sum_{k\in{\cal I}_j^{in}}\min\left\{\alpha_{lk}^jD_{k}(\rho_{k}(L_{k},.)),\,\beta_{kl}^j S_{l}(\rho_{l}(0,.))\right\} \quad \mbox{ in }(0,T), \label{MT_CCj_out}\\
f_z(\rho_z(L_z,.))=\min\{f_z^{out},D_z(\rho_z(L_z,.))\} \quad \mbox{ in }(0,T), \label{MT_OutC}\\
f_y(\rho_y(0,.))=\min\{D_y^{in}(q_y,.),S_y(\rho_y(0,.))\} \quad \mbox{ in }(0,T), \label{MT_InC}
\end{align}
\begin{equation}\label{PVI_queue}
\hspace*{3cm}\left.\begin{array}{l}
\displaystyle\frac{dq_y}{dt}=f_y^{in}-f_y(\rho_y(0,.)) \quad \mbox{ in }(0,T),\\
q_y(0)=q_y^0,
\end{array}\right\}
\end{equation}\label{MT}
\end{subequations}

\noindent where the terms $D_i$ and $S_i$ represent the demand and supply functions respectively \cite[see][]{alvarez-vazquez2017, vazquez-mendez2019}, the term $D_y^{in}(q_y,t)$ represents the demand of queue $q_y$ at time $t$  \cite[cf.][]{vazquez-mendez2019}, and values $\rho_i^0$ and $q_y^0 \geq 0$ are, respectively, the initial density at road $A_i$ and the initial queue length downstream avenue $A_y$. \\  

Moreover, 
\begin{itemize}
\item the parameters $\alpha^j_{lk}$ stand for drivers' preferences when arriving at a junction, i.e., $\alpha^j_{lk}$ represents the rate of drivers that, reaching intersection $j$ coming from road $A_k$, will take the outgoing road $A_l$. Thus, these compatibility constraint need to be verified:
\begin{equation}\label{constraint1}
0\leq\alpha^j_{lk}\leq 1 \quad \mbox{ and } \quad \sum_{l\in{\cal I}_j^{out}}\alpha^j_{lk}=1.
\end{equation}
\item the parameters $\beta^j_{kl}$ stand for ingoing capacities at outgoing roads, i.e., $\beta^j_{kl}$ represents the rate of vehicles that, arriving at junction $j$ for road $A_k$, can enter the outgoing road $A_l$. As above, these parameters should satisfy:
\begin{equation}\label{constraint2}
0\leq\beta^j_{kl}\leq 1 \quad \mbox{ and } \quad \sum_{k\in{\cal I}_j^{in}}\beta^j_{kl}=1.
\end{equation}
\end{itemize}
Finally, it is worthwhile recalling here the fundamental role of coupling conditions (\ref{MT_CCj_in}) and (\ref{MT_CCj_out}) in order to guarantee the conservation of the number of cars at intersections.

\subsection{Modelling atmospheric pollution}

Traffic realted air pollution is simulated here by a mathematical model similar to the one proposed in \cite{alvarez-vazquez2018}, already used in \cite{vazquez-mendez2019} and whose uniqueness of solution was argumented in \cite{ca2,martinez2017}. Due to its main role, we focus our interests only in pollution related to nitrogen oxides (NO$_x$), but many other kinds of pollution -like carbon monoxide (CO), sulphur oxides (SO$_x$), total hydrocarbons (THC), etc.- could be also included. So, the NO$_x$ concentration $\phi(x,t)$ $[kg/km^2]$ corresponding to each point $x\in\Omega$ and each time $t\in[0,T]$, can be obtained by solving the following initial/boundary value problem:
\begin{subequations}
\begin{eqnarray}
\displaystyle \frac{\partial \phi}{\partial t} + \mathbf{v}\cdot\nabla\phi - \nabla\cdot (\mu\nabla\phi)  + \kappa\phi = \sum_{i=1}^{N_R}\xi_{A_i}  \quad \mbox{ in }\Omega\times (0,T), \label{polucion} \\
\displaystyle \phi(.,0) = \phi^0  \quad \mbox{ in }\Omega , \label{polucion_ci} \\
\displaystyle \mu\frac{\partial\phi}{\partial n} - \phi\,\mathbf{v}\cdot\mathbf{n} =\sum_{y\in{\cal I}^{in}} \lambda_yq_y\delta_{\sigma_y(0)} \quad \mbox{ on }S^- , \label{polucion_cc1} \\
\displaystyle \mu\frac{\partial\phi}{\partial n} = 0  \quad  \mbox{ on }S^+, \label{polucion_cc2}
\end{eqnarray}\label{MP}
\end{subequations}

\noindent where the field $\mathbf{v}(x,t)$ $[km/h]$ denotes wind velocity, the function $\phi^0$ is the given initial NO$_x$ concentration, the coefficients $\mu (x,t)$ $[km^2/h]$ and $\kappa(x,t)$ $[h^{-1}]$ are, respectively, the NO$_x$ molecular diffusion and the NO$_x$ extinction rate, the terms $\lambda_yq_y\delta_{\sigma_y(0)}$ represent pollution sources due to queues entering by the inflow boundary (cf. \cite{vazquez-mendez2019} for further details), vector $\mathbf{n}$ denotes the unit outward normal vector to the boundary $\partial \Omega = S^-\cup S^+$, split into the outflow boundary
$S^+ = \{(x,t)\in \partial\Omega\times (0,T) \mbox{ such that } \mathbf{v}\cdot\mathbf{n} \geq 0 \}$ and the inflow boundary $S^- = \{(x,t)\in \partial\Omega\times (0,T) \mbox{ such that } \mathbf{v}\cdot\mathbf{n} < 0 \}$. Moreover, the terms $\xi_{A_i}$ $[kg/km^2/h]$ represent pollution sources due to vehicular traffic on roads $A_i$, and are given by means a Radon measure: For each $t\in[0,T]$, the distribution $\xi_{A_i}(t):\mathcal{C}(\overline\Omega)\longrightarrow \mathbb{R}$ is defined by:
$$
\langle \xi_{A_i}(t),v \rangle =\int_{0}^{L_i} \left(\gamma_i f_i(\rho_i(s,t))+\eta_i\rho_i(s,t)\right)\,v(\sigma_i(s))\, ds, \quad \forall v\in \mathcal{C}(\overline\Omega),
$$
where $\sigma_i$ is the parametrization of avenue $A_i$, density $\rho_i$ is given by the traffic model (\ref{MT}), and parameters $\gamma_i$ and $\eta_i$ are weights associated to pollution rates.

\section{A bi-level non-cooperative optimal control problem}

In previous approaches to traffic management in a road network, standard objectives were usually related only to traffic problems, such as travel time or congestions. Nevertheless, present-day difficulties with air pollution in the surroundings of big metropolises have turned the mitigation of this phenomenon into another major aim in the optimal management of urban road networks.

Here, two different objectives, one of each type, will be considered in a simultaneous way. With respect to optimizing the traffic flow, it is important to minimize the total travel time and to maximize the outflow of the network. In the spirit, for instance, of \cite{goatin2016,vazquez-mendez2019}, the following functional $J_T$ should be minimized:
\begin{equation}\label{coste_trafico}
\hspace*{-0.25cm}J_T=\int_0^T\hspace*{-0.15cm}\left(\sum_{y\in{\cal I}^{in}}\epsilon^q_y q_y(t)+\sum_{i=1}^{N_R}\epsilon_i\int_{0}^{L_i}\rho_i(s,t)\,ds -\hspace*{-0.2cm}\sum_{z\in{\cal I}^{out}}\hspace*{-0.1cm}\epsilon^{out}_z f_z(\rho_z(L_z,t))\hspace*{-0.1cm}\right) dt,
\end{equation}
where $\epsilon^q_y,\, \epsilon_i,\,\epsilon^{out}_z\geq 0$ are weight parameters to be chosen by the decision makers according to their preferences.

Regarding air pollution, it is essential to keep mean concentration of NO$_x$ as low as possible, i.e., we are involved in minimizing the cost functional $J_P$ giving the mean pollution concentration:
\begin{equation}\label{coste_polucion}
J_P=\frac{1}{T\, |\Omega|}\int_0^T \int_\Omega \phi(x,t)\,dx\,dt,
\end{equation}
where $|\Omega|$ denotes the usual Euclidean measure of set $\Omega$. (We must remark here that the averaged value could be taken in any sensitive region $D \subset \Omega$ and in any time subinterval of $[0,T]$ but, for the sake of simplicity, we have chosen here the full domains).

For the {\it controls} (that is, the design variables that can be managed within the network), several different choices have been investigated in previous studies: incoming fluxes \cite{goatin2016}, drivers' preferences \cite{gugat}, network expansions \cite{alvarez-vazquez2018}, etc. However, we will center our attention at the optimal management of the network intersections, attempting to obtain those parameters $\alpha^j_{lk}$, $\beta^j_{kl} $ that are the most satisfactory for our global aims. 

Supposing that the parameters $ \alpha_{lk}^j $ (drivers' preferences) change when the input/output ratios $ \beta_{kl}^j $ are modified at the intersections, we will assume that the set of drivers always try to minimize the functional $ J_T $, while the leader organization managing the whole network intends to choose the ratios to try to minimize atmospheric contamination.

Following this reasoning, we face up to a bi-level problem. In a firts level, we have the {\it follower problem}: For a given $ \mathbf{\beta}=(\beta_{kl}^j), \ j=1,\ldots,N_J,\ k\in{\cal I}_j^{in},\ l\in{\cal I}_j^{out}$, verifying (\ref{constraint2}), solve:
\begin{equation}\label{P2f}
\begin{array}{l}
\displaystyle  \min J_T(\mathbf{\alpha},\mathbf{\beta} )\\
\mbox{subject to (\ref{constraint1})}
\end{array}
\end{equation}
with $ \mathbf{\alpha}=(\alpha_{lk}^j), \ j=1,\ldots,N_J,\ k\in{\cal I}_j^{in},\ l\in{\cal I}_j^{out}. $

Then in a second level, the {\it leader problem} reads as:
\begin{equation}\label{P2l}
\begin{array}{l}
\displaystyle  \min J_P(\mathbf{\alpha}_{\mathbf{\beta}},\mathbf{\beta})\\
\mbox{subject to (\ref{constraint2})}
\end{array}
\end{equation}
where $\mathbf{\alpha}_{\mathbf{\beta}}$ is the optimal solution of the follower problem (\ref{P2f}) for given data $\beta$.

In this approach, our main objective relies in computing a Stackelberg strategy for the bi-level problem (\ref{P2f})-(\ref{P2l}), in the sense below classical definition:

\begin{definition}\label{DefStackelberg}
A pair $(\mathbf{\alpha}^*,\mathbf{\beta}^*)$ is said to be a Stackelberg strategy, solution of the bi-level problem (\ref{P2f})-(\ref{P2l}), if it verifies that:
\begin{enumerate}
\item $\mathbf{\alpha}^*$ is the best reaction of the follower to the leader choice $\mathbf{\beta}^*$, i.e., $\mathbf{\alpha}^*$ is the solution of the follower problem (\ref{P2f}) for given data $\mathbf{\beta}^*$ (in other words, $\mathbf{\alpha}^*=\mathbf{\alpha}_{\mathbf{\beta}^*}$).
\item $\mathbf{\beta}^*$ is the best option of the leader, i.e., $\mathbf{\beta}^*$ is the optimal solution of the leader problem (\ref{P2l}).
\end{enumerate} 
\end{definition}

We must remark here that, by using adjoint techniques \cite{marchuk}, the functional $J_P(\mathbf{\alpha},\mathbf{\beta} )$ can be rewritten in the more useful alternative form (see full details in Theorem 3.1 of \cite{alvarez-vazquez2018}):
\begin{equation}\label{JP_adjoint}
\begin{split}
J_P = \sum^{N_R}_{i=1} \int^T_0 \int^{L_i}_{0} (\gamma_i f_i(\rho_i(s,t)) + \eta_i\rho_i(s,t) ) \, g(\sigma_i(s),t) \, ds \, dt \\ 
+ \sum_{y\in\mathcal{I}^{in}} \int^T_0 \lambda_y q_y(t) \, g(\sigma_y(0),t) \, \chi_{S^-}(\sigma_y(0),t) \, dt + \int_\Omega \phi^0(x) \, g(x,0) \, dx,
\end{split}
\end{equation}
where $\chi_{S^-}$ is the characteristic function of the inflow boundary $S^-$, $ \gamma_i f_i(\rho_i(s,t))$ $ + \eta_i\rho_i(s,t) $ represents the pollutant emissions on the road network, and $g(\sigma_i(s),t)$ is the evaluation on the road network of the so-called {\it adjoint state} $g(x,t)$, the solution of the following final/boundary value problem: 
\begin{subequations}\label{Adjoint_model}
\begin{align}
-\frac{\partial g}{\partial t} - \mathbf{v}\cdot\nabla g - \nabla\cdot(\mu\nabla g) + \kappa g = \frac{1}{T \, |\Omega|} \quad \mbox { in } \Omega\times(0,T), \\
g(x,T) = 0 \quad \mbox{ in } \Omega , \\
\mu\frac{\partial g}{\partial n} = 0 \quad \mbox{ on } S^- , \\
\mu\frac{\partial g}{\partial n} +  g\mathbf{v}\cdot\mathbf{n} = 0 \quad \mbox{ on } S^+ .
\end{align}
\end{subequations}   

This alternative formulation of functional (\ref{JP_adjoint}) depends straightforwardly on traffic density and flow rate and, consequently, on the controls $(\alpha,\beta)$. However, this control dependency is implicit (as can be seen in conditions (\ref{MT_CCj_in})-(\ref{MT_CCj_out})) making hard to get an explicit expression of the derivative of the functionals (\ref{coste_trafico}) and (\ref{JP_adjoint}) with respect to the controls. This fact will be a decisive issue in the choice of the methods for solving the bi-level multi-objective control problem.

Finally, it is important emphasizing here that the adjoint state $g(x,t)$ is the unique solution of the adjoint equation \cite[cf.][]{lsu,alvarez-vazquez2018}, being independent of the traffic variable $\rho(x,t)$. Consequently, the adjoint state does not depend on the traffic model. So, it can be computed separately, and adjoint problem (\ref{Adjoint_model}) only needs to be solved once in a preliminary step. 

\section{Numerical solution of the bi-level multi-objective control problem}

The bi-level problem (\ref{P2f})-(\ref{P2l}) is generally non-convex. Therefore, many local solutions are expected. Moreover, effective expressions for the gradients of objective functionals $J_T$ and $J_P$ with respect to the controls $(\alpha,\beta)$ are hard to compute (leading to the only reasonable option involving its numerical approximation). So, free-derivative optimization methods or methods using numerical approximation of gradients will be the natural and efficient choice in order to solve the bi-level problem. 

\subsection{Discretization of cost functionals $J_T$ and $J_P$}

With independence on the method chosen to solve the bi-level problem, its efficiency relies on a good discretization and evaluation of the cost functionals $J_T$ and $J_P$. So, as a previous step, we show how this can be performed in a suitable way (following the method already introduced in \cite{vazquez-mendez2019}).
 
We choose the following space-time discretization: For each road $A_i$, the
parametrization interval $I_i = [0,L_i]$ is split into $M_i$ cells $I_{i,h} = [s_{i,h-\frac{1}{2}}, s_{i,h+\frac{1}{2}}], \ h = 1, \ldots , M_i$, of length $\Delta s_i > 0$, where $s_{i,h} = (s_{i,h-\frac{1}{2}} + s_{i,h+\frac{1}{2}})/2$ represents the midpoint of each cell. On the other part, the time interval $[0, T ]$ is also split into $N \in \mathbb{N}$ subintervals of length $\Delta t = T /N$, defining in this way the discrete times $t^n = n\Delta t, \ n = 0, \ldots , N$, are defined. Using this discretization the system (\ref{MT}) can be solved addressing the functional $J_T$ with quadrature rules \cite{alvarez-vazquez2015, alvarez-vazquez2017}. In particular, given the discrete density $\rho^n_{i,h}$ and queue $q^n_i$, for $n=0,\ldots,N,\ i=1,\ldots,N_R, \ h = 1,\ldots,M_i$, we evaluate the following full-discrete integral: 
\begin{equation}
J^\Delta_T = \Delta t \sum^N_{n=0}\left( \sum_{y\in\mathcal{I}^{in}} \epsilon^q_y q^n_y +
\sum^{N_R}_{i=1}\epsilon_i \Delta s_i \sum^{M_i}_{h=1}\rho^n_{i,h} - \sum_{z\in\mathcal{I}^{out}} \epsilon^{out}_z f_z(\rho^n_{z,M_z}) \right).
\end{equation}
On the other part, let us consider a polygonal approximation $\Omega_{h}$ of $\Omega$, with an admissible triangulation $\tau_{h}$, where vertices $\{x_{j}, \ j=1,\ldots,N_{v}\}$ satisfy that all the vertices on the boundary $\partial\Omega_{h}$ remain on the boundary $\partial\Omega$, that is $\sigma_y(0),\,\sigma_z(L_z)\in\partial\Omega_h,$ for all $y\in{\cal I}^{in},\, z\in{\cal I}^{out}$, and that, for $n=0,\ldots,N-1$, each edge of $\partial\Omega_h$ lies either in $({S_h^n})^-=\{x\in\partial\Omega_h \ : \ \mathbf{v}\cdot\mathbf{n} < 0\}$ or in $({S_h^n})^+=\{x\in\partial\Omega_h \ : \ \mathbf{v}\cdot\mathbf{n}\geq 0\}$.\\
Then, the adjoint model can be solved numerically on the domain $\Omega_h$, and for the discrete times $\{t^n\}^N_{n=0}$ we get the discrete adjoint values $\{\{g^n_{h,k}\}^{n_v}_{k=0} \}^N_{n=0}$ (see \cite{alvarez-vazquez2015} and Algorithm 3 of \cite{vazquez-mendez2019}). Once this is done, it is possible to evaluate the adjoint state at roads' nodes getting $\{g^n_h(\sigma_i(s_{i,h}))\}^{N_R}_{i=1}$ by triangular interpolation. Thus, the leader functional $J_P$ can be now addressed by quadrature rules: Given the discrete functions $\rho^n_{i,h}$, $f_i(\rho^n_{i,h})$, $g^n_{h,k}$ and $g^n_h(\sigma_i(s_{i,h}))$, for $n=1,\ldots,N, \ i=1,\ldots,N_R,\  h=1,\ldots,M_i$, we evaluate the following full-discrete integral:  
\begin{equation}
\begin{split}
J^\Delta_P = \Delta t\sum^N_{n=1}\sum^{N_R}_{i=1}\sum^{M_i}_{h=1} \Delta s_i (\gamma_i f_i(\rho^n_{i,h}) +\eta_i\rho^n_{i,h} ) g^n_h(\sigma_i(s_{i,h})) \Vert \sigma'_i(s_{i,h})\Vert \\
+ \sum^N_{n=1} \sum_{\substack{y\in\mathcal{I}^{in}\\
                            \sigma_y(0)\in (S^n_h)^-
                   }} \lambda_y q^n_y g^n_h(\sigma_y(s_{y,1})) 
                   + \frac{1}{3} \sum_{\mathcal{T}\in \tau_h} |\mathcal{T}|\sum_{x_j\in\mathcal{T}} \Phi^0(x_j)g^0_j
\end{split}                                       
\end{equation} 
It is worthwhile remarking here that each evaluation of the discrete functionals $J^\Delta_T$ and $J^\Delta_P$ for a pair $(\alpha,\beta)$ requires the solution of the LWR traffic model (\ref{MT}) (for example, by Algorithm 1 of \cite{vazquez-mendez2019}), the computation of the adjoint state (Algorithm 3 of \cite{vazquez-mendez2019}), and also efficient algorithms to evaluate the objective functionals (see, for instance, Algorithms 2 and 4 of \cite{vazquez-mendez2019}). 

\subsection{Solving the follower problem}

Given the discrete cost functional $J^\Delta_T$, the follower problem (\ref{P2f}) will be solved by combining an interior-point method and a genetic algorithm. Both algorithms are implemented respectively by the solvers \texttt{fmincon} and \texttt{ga} from the Optimization Toolbox of Matlab R2017a, being important the following issues: The solver \texttt{ga} includes a \texttt{hybrid} option that allows combining it with other Matlab optimization solvers, can be executed in parallel, and uses the three basic {\it probabilistic rules} of the natural selection: elite, crossover and mutation to generate the next generation (cf. Algorithm 1). On the other part, the solver \texttt{fmincon} can approximate the cost functional gradient in case of not availability (like our case), and can be also executed in parallel. Moreover, the direction-search of solver \texttt{ga} presents a large set of probabilities, provided by the population diversity; in contrast, in the solver \texttt{fmincon} the direction-search is given by a line-search and trust-region criterion which depends of the direct-step or conjugate-gradient step. Then, with the aim of providing  \texttt{fmincon} with a diversity similar to \texttt{ga}, in this work a multi-start execution of \texttt{fmincon} was carried out (see full details in Algorithm 2).

\begin{algorithm}
\SetAlgoLined
\KwData{Initial vectors population $\tilde\alpha^0=\{\alpha^{0,n} \}^N_{n=1}$, fixed vector $\beta$, and tolerance $tol$}
\KwResult{ Optimal vector of preferences $\alpha^*$, and optimal functional value $J^\Delta_T(\alpha^*,\beta)$ }
\Begin{
 \texttt{set} $k=0$\;
 \While{$Error>tol$}{
  \For{$n = 1,...,N$ }{ 
      \texttt{Compute} $J^\Delta_T(\alpha^{k,n}, \beta)$ by \texttt{Algorithm 2} of \cite{vazquez-mendez2019} 
      }
      \texttt{Generate} the new population $\tilde\alpha^{k+1}$ by natural selection: \texttt{Elite, Crossover} and \texttt{Mutation}\; 
      \texttt{Compute} $Error$ and \texttt{set} $k = k+1$\;
  }
  \texttt{set} $\bar\alpha$ as preference corresponding to the mean of functional value set $\{ J_T(\alpha^{k+1,n}, \beta) \}_{n=1}^N$\;
  \texttt{switch} to \texttt{fmincon}\;
  \texttt{input} $\bar\alpha$ and \texttt{compute} $\min J^\Delta_T(\alpha,\beta)$ by \texttt{fmincon}, and \texttt{get} the optimal $\alpha^*$ and the optimal functional value  $J^\Delta_T(\alpha^*,\beta)$\; 
  }
  \caption{\texttt{ga} algorithm with hybrid option }
  \end{algorithm}
  
\begin{algorithm}
\KwData{Multi-initial vectors set $\alpha = \{\alpha^i_0\}^N_{i=1}$, fixed vector $\beta$, and tolerance $tol$}
\KwResult{Optimal vector of preferences $\alpha^*$, and optimal functional value $J^\Delta_T(\alpha^*,\beta)$ }
\SetAlgoLined
\Begin{
\For{$i=1,...,N$}{
    \texttt{Input} $\alpha^i_0$ and \texttt{compute} $\min J^\Delta_T(\alpha,\beta)$ by \texttt{fmincon} and get the optimal $\alpha^{*,i}$   
    } 
    \texttt{set} $J^\Delta_T (\alpha^*,\beta) = \min \{J^\Delta_T(\alpha^i,\beta)\}^N_{i}$\;
    \texttt{set} $\alpha^*$ as the best from $\{\alpha^{*,i}\}^N_{i=1}$\;  
 }
 \caption{\texttt{fmincon multi-start} algorithm}           
\end{algorithm}

\subsection{Solving the leader problem}

To compute the Stackelberg solution $(\alpha^*,\beta^*)$ in the sense of Definition 1, a combination between the solvers \texttt{ga-hybrid} (Algorithm 1) and \texttt{fmincon} multi-start (Algorithm 2) will be used, with the aim of addressing the complexity of the leader problem. Thus, the diversity at searching directions provided by \texttt{ga} could be used for identifying a feasible initial value for the \texttt{fmincon}, getting in this way a high quality Stackelberg solution. Nevertheless, since this hybrid method only gives $\beta^*$ as output, the best follower response to the leader $\alpha^*$ needs to be calculated and saved in the last evaluation of the leader functional at the \texttt{fmincon} stage of the hybrid solver. All the details of this process are shown in Algorithm 3, where the use of the adjoint state provides an important saving in the total computational cost.

\begin{algorithm}
\KwData{Initial vectors population $\tilde\beta^0=\{\beta^{0,n}\}^N_{n=1}$}
\KwResult{Stackelberg solution $(\alpha^*,\beta^*)$, and optimal functional values $J^\Delta_T(\alpha^*, \beta^{*}), J^\Delta_P(\alpha^*, \beta^{*})$}
\SetAlgoLined
\Begin{
 \texttt{set} $k=0$ and \texttt{compute} the adjoint $g^n_{h,k}$ by \texttt{Algorithm 3} of \cite{vazquez-mendez2019}\;
 \While{$Error>tol$}{
  \For{$n = 1,...,N$ }{ 
   \texttt{set} randomly $\alpha = \{\alpha^i_0\}^M_{i=1}$  as multi-initial vectors set\;
   \texttt{Input} $\alpha$ and \texttt{compute} $\min J^\Delta_T(\alpha, \beta^{k,n})$ by \texttt{Algorithm 2}, and get $\alpha_{\beta^{k,n}}$\; 
   \texttt{Compute} $J^\Delta_P(\alpha_{\beta^{k,n}},\beta^{k,n})$ by \texttt{Algorithm 4} of \cite{vazquez-mendez2019}  \;
  }
  \texttt{Generate} the next population $\tilde\beta^{k+1}$ by natural selection: \texttt{Elite, Crossover and Mutation}\; 
  \texttt{Compute} $Error$ and \texttt{set} $k = k+1$\;
  }
  \texttt{set} $\bar\beta$ as the vector of restrictions corresponding to the mean of the functional values set $\{J_P(\alpha, \beta^{k+1,n})\}^N_{n=1}$\;
  \texttt{switch} to \texttt{fmincon}\;
  \texttt{input} $\bar\beta$ and \texttt{compute} $\min J^\Delta_P(\alpha,\beta)$ by \texttt{fmincon}, and \texttt{get} $\beta^*$\;
  from the last evaluation of $J^\Delta_P(\alpha,\beta)$ \texttt{get} $\alpha^*$\; 
  \texttt{Compute} $J^\Delta_T(\alpha^*, \beta^{*})$ and $J^\Delta_P(\alpha^*, \beta^{*})$
 }  
 \caption{Stackelberg algorithm}
\end{algorithm}

\section{Numerical experiences}

We present and analyse here several computational results obtained in a real-world scenario in Mexico, set in the Guadalajara Metropolitan Area (GMA). Given the previous experiences over the same domain developed by the authors in \cite{alvarez-vazquez2017, alvarez-vazquez2018, vazquez-mendez2019} we only present here essential data and assumptions in a summarized way.

\subsection{Initial/boundary conditions and models parameters}

The road network analyzed here is composed by $N_R=15$ avenues and $N_J = 9$ junctions, where all avenues have only one lane and its theoretical flow is given by the static relation defined in \cite{alvarez-vazquez2018}. As boundary conditions for the traffic model (\ref{MT}), we consider equal downstream road capacities for the 3 incoming avenues ($C^{in}_i = 2.013\,10^3, \ i=1,2,10$), with equal sinusoidal desired inflow rate, and also equal maximum outflow rates for the 3 outgoing roads ($f^{k}_z = 2.013\,10^3, \ k=13,14,15$). As initial conditions, null traffic ($\rho^0_{i,s} = 0, \ i=1,...,N_R$) and queues ($q^0_y = 0, \ y = 1,2,10$) were assumed. Also, we have considered the following weights parameters: For the road densities $\epsilon_i = 0.7, \ i=1,2,3$ and $\epsilon_i = 0.5$ for the rest of avenues, for queue lengths $\epsilon^q_i = 0.45, \ i = 1,2$ and $\epsilon^q_{10} = 0.1$, and for outflow rates $\epsilon^{out}_z = 0.5,\ z = 13,14,15$. This combination of weights values translates the intention of increasing the follower cost by an excess of traffic densities and/or queue lengths. 

\begin{figure}
\begin{center}
\includegraphics[scale=1]{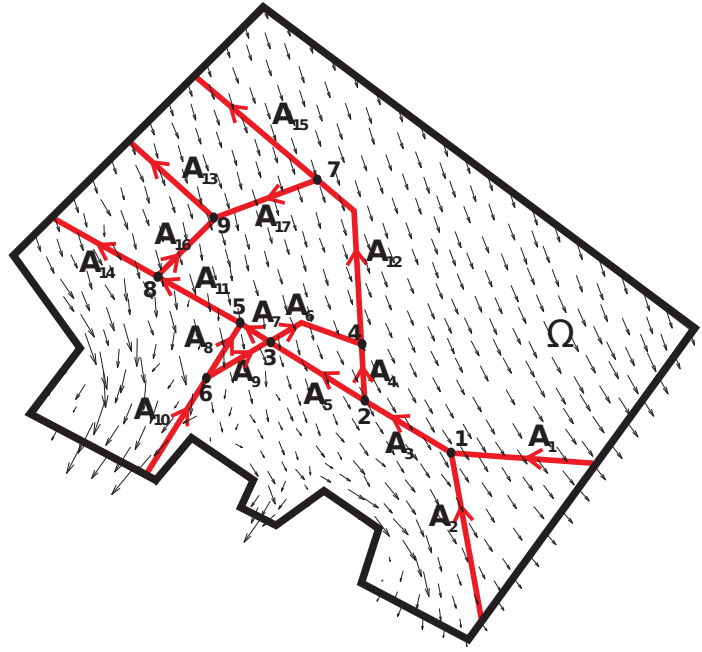}
\end{center}
\caption{The polygonal domain $\Omega $ made it over a satellite image of the GMA (no displayed here). The vectors depicting the wind field are drawn in black and the straight lines depicting the road network are drawn in red.}\label{Fig4}
\end{figure}

With respect to the pollution model (\ref{MP}) and its corresponding adjoint state (\ref{Adjoint_model}), typical values for NO$_x$ ($\mu = 3.5\, 10^{-8}\, km^2/h, \kappa = 0.6 10^{-2} h^{-1}, \gamma_i = 10^6\, kg/$ $number\, of\, cars/km, \eta_i = 3.16 10^{-5}$ $kg/number\, of\, cars/h$) have been taken. We also consider null pollution at initial time, not external pollution sources and, due to the particular wind direction at $S^+$ (see Fig. \ref{Fig4}), parameters $\lambda_y = 0$.

Regarding the discretization, we chose a time step of $\Delta t = 4\, 10^{-3}$ (measured in hours) and, for each road $A_i$, its spatial domain $I_i$ has been divided into cells large enough to guarantee the classical CFL condition ($\Delta s_i \in (0.2, 0.21)$). Finally the polygonal domain $\Omega_h \subset R^2$ has been discretized with a triangulation of 898 triangles and 491 vertices, satisfying standard regularity hypothesis, in order to guarantee the numerical method convergence \cite{GMSH}.\\
 With respect to the minimization algorithms, we consider an initial population of $50$ individuals for the \texttt{ga} and \texttt{ga-hybrid} solvers, and a set of 5 vectors as initial input for the \texttt{fmincon-multi-start} solver. All solvers have been executed in parallel in an AMD Threadripped 1920X CPU at 3.8 GHz with 12 cores and 24 threads desktop, 32 GB RAM, and Linux Mint OS.    


\subsection{Assessment experiences for solving the follower problem}

In the follower problem (\ref{P2f}) the restrictions vector $\beta= \{\beta^{j}_{kl}\}$ must be fixed and, in this particular experience, it will be chosen in such a manner that we can predict the preferences vector $\alpha = \{\alpha^j_{lk}\}$ and, consequently, can confirm the reliability of our approach. Then, in this spirit, two different representative cases are shown.

\underline{Case 1:} In this first case, the restrictions $\beta^{1,j}_{kl}$ are taken such that the avenues $A_5, A_6$ and $A_7$ remain blocked at intersections $j = 3, 4$ and $5$ (that is, $\beta^{1,3}_{5\,6} = \beta^{1,3}_{5\,7} = \beta^{1,4}_{6\,12} = \beta^{1,5}_{7\,11} = 0$). Then, it is expected that the drivers coming from $A_3$ and $A_{10}$ take the avenues $A_4$ and $A_{8}$, respectively, avoiding the block imposed by the leader. In this case it is also expected that drivers take their respective outways by avenues $A_{12}$ and $A_{11}$. Thus, vectors $\beta^{1,j}_{kl}$ were fixed satisfying above constraints (see Table \ref{TablaFollower}) and then the follower problem was solved using  \texttt{ga-hybrid}, \texttt{ga} and \texttt{fmincon-multi-start} routines. These solvers were addressed with a tolerance of $10^{-4}$ (although in the hybrid case an additional tolerance of $10^{-10}$ was imposed), and the corresponding resulting preferences, respectively denoted by $\alpha^{hyb,j}_{lk}, \alpha^{ga,j}_{lk}$ and $\alpha^{fmin,j}_{lk}$, are shown at Table \ref{TablaFollower}.

\begin{figure}
\centering
\subfloat[Case 1.]{%
\resizebox*{5cm}{!}{\includegraphics{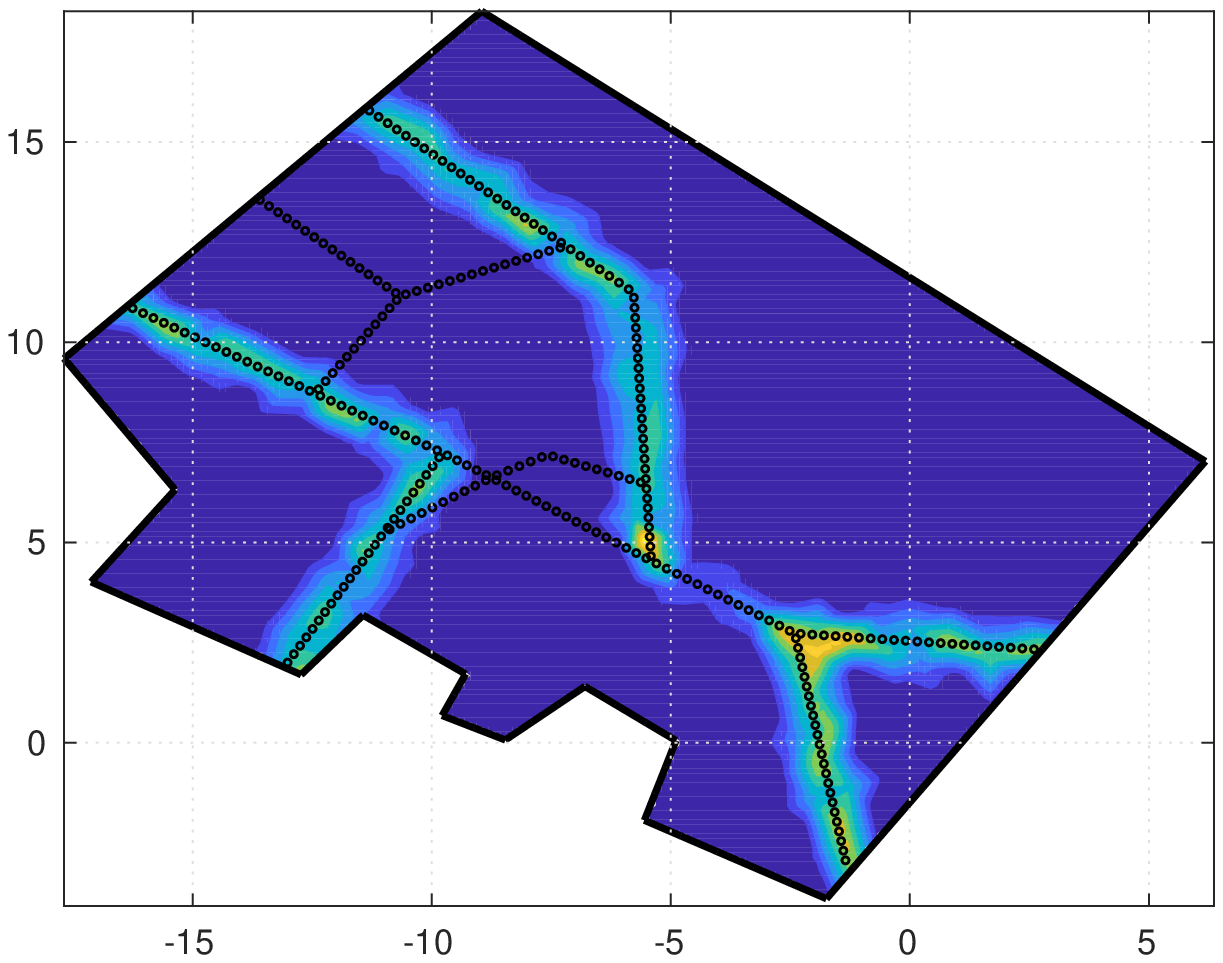}{\label{Caso1}}}}\hspace{5pt}
\subfloat[Case 2.]{%
\resizebox*{5cm}{!}{\includegraphics{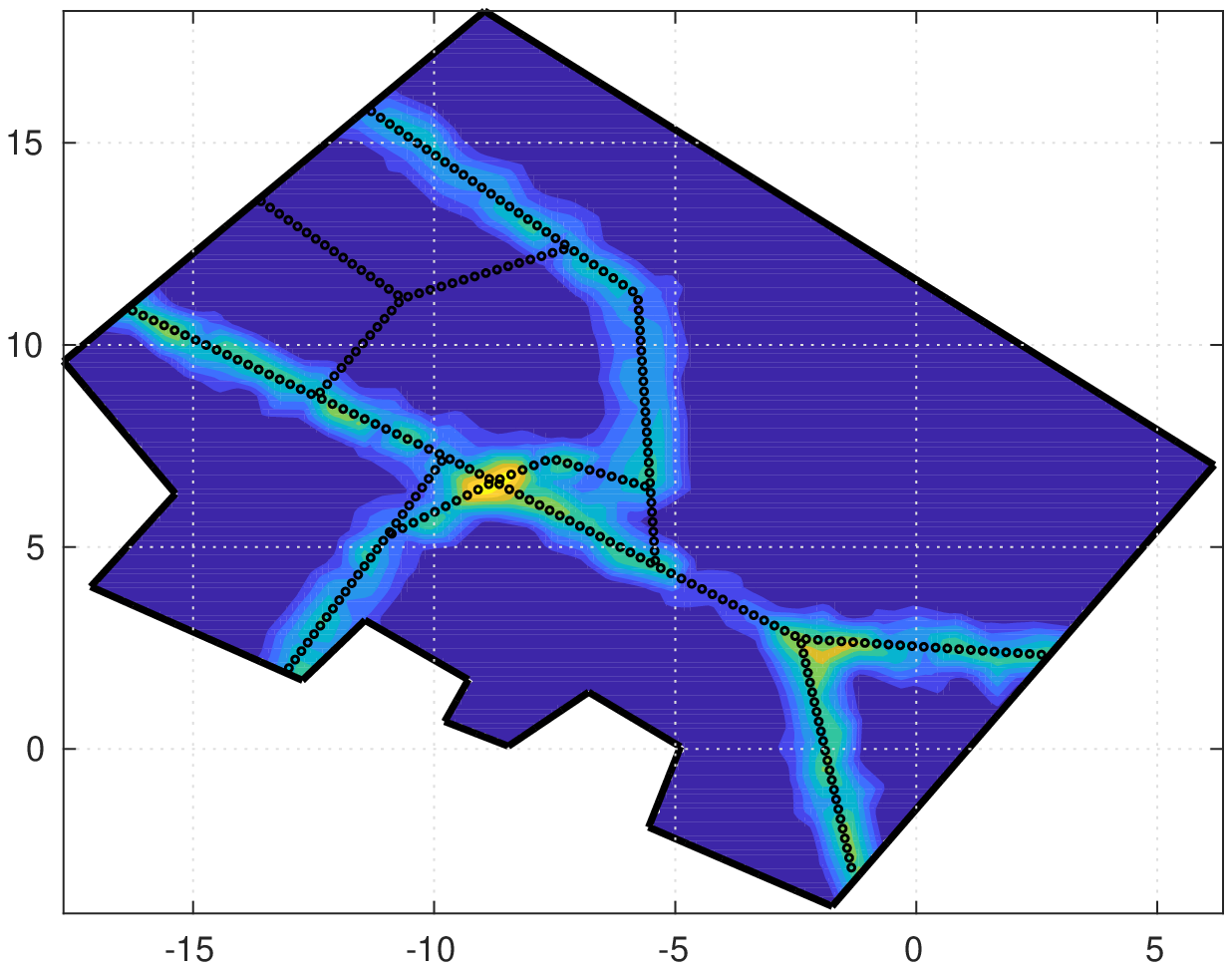}{\label{Caso2}}}}
\caption{Level curves of the mean concentrations of pollution for a simulation period of 24 hours with zero wind.}
\label{Isolineas_NoOptimo1}
\end{figure}

Here, the three solutions present in common the preferences $\alpha^2_{4\,3} = 1, \alpha^6_{8\,10} = 1$, which imply that all drivers from $A_3$ and $A_{10}$ take the avenues $A_4$ and $A_{8}$ at junctions $j=2$ and $j=6$, avoiding the blocked avenues $A_5, A_6$ and $A_7$. To confirm this, in Figure \ref{Caso1} the isolines of pollution concentration (with zero wind $ \mathbf{v}=0$) are depicted, showing in a clear way in which avenues the pollutant emissions are present (and consequently the traffic flow is high). Therefore, the predicted behavior of drivers is fulfilled for this first case.

\begin{table}
\tbl{Data for the two analyzed cases of the follower problem (\ref{P2f}): The fixed leader restrictions $\beta^{1,j}_{kl}$ and $\beta^{2,j}_{kl}$ correspond to Case 1 and Case 2, respectively. For each case, the optimal preferences $\alpha^{solver,j}_{lk}$ for the three solvers (\texttt{fmincon, hybrid, ga}) are also displayed.}
{\begin{tabular}{|c|ccccccccc|}
   \hline
                           & j=1        & j=2            & j=3         & j=4       &  j=5      & j=6              & j=7         & j=8         & j=9 \\
   \hline
 $\mathcal{I}_j^{in}(k)$   & $\{1,2\}$  & $\{3\}$        & $\{5,9\}$   & $\{4,6\}$ & $\{7,8\}$ & $\{10\}$         & $\{12\}$    & $\{11\}$    & $\{16,17\}$\\
 $\mathcal{I}_j^{out}(l)$  & $\{3 \}$   & $\{4,5\}$      & $\{6,7 \}$  & $\{12\}$  & $\{11\}$  & $\{8,9\}$        & $\{15,17\}$ & $\{14,16\}$ & $\{13\}$   \\
    \hline
 $\beta^{1,j}_{kl}$        & $\begin{bmatrix}
                             0.5\\
                             0.5 
                        \end{bmatrix}$  & $(1,1)$       & $\begin{bmatrix}
                                                          0.00 & 0.00\\
                                                          1.00 & 1.00
                                                         \end{bmatrix}$ & $\begin{bmatrix}
                                                                           1.00 \\
                                                                           0.00
                                                                           \end{bmatrix}$ & $\begin{bmatrix}
                                                                                       0.00 \\
                                                                                       1.00
                                                                                       \end{bmatrix}$ & $(1,1)$ & $(1,1)$ & $(1,1)$ & $\begin{bmatrix}
                                                                                                                                        0.5\\
                                                                                                                                        0.5
                                                                                                                                       \end{bmatrix}$\\                                                                                                                                                                                                                                                                            
$\alpha^{hyb,j}_{lk}$       & $(1,1)$ & $\begin{bmatrix}
                                         1.00 \\
                                         0.00
                                        \end{bmatrix}$ & $\begin{bmatrix}
                                                       0.40 & 0.16\\
                                                       0.60 & 0.84
                                                      \end{bmatrix}$ & $(1,1)$  & $(1,1)$    & $\begin{bmatrix}
                                                                                               1.00\\
                                                                                               0.00
                                                                                                \end{bmatrix}$ & $\begin{bmatrix}
                                                                                                                 1.00\\
                                                                                                                 0.00
                                                                                                                 \end{bmatrix}$ & $\begin{bmatrix}
                                                                                                                                  1.00\\
                                                                                                                                  0.00 
                                                                                                                                  \end{bmatrix}$ & $(1,1)$ \\
$\alpha^{ga,j}_{lk}$       & $(1,1)$ & $\begin{bmatrix}
                                          1.00 \\
                                          0.00
                                          \end{bmatrix}$ & $\begin{bmatrix}
                                                           0.38 & 0.02\\
                                                           0.62 & 0.98
                                                      \end{bmatrix}$ & $(1,1)$  & $(1,1)$    & $\begin{bmatrix}
                                                                                               1.00\\
                                                                                               0.00
                                                                                                \end{bmatrix}$ & $\begin{bmatrix}
                                                                                                                 1.00\\
                                                                                                                 0.00
                                                                                                                 \end{bmatrix}$ & $\begin{bmatrix}
                                                                                                                                  1.00\\
                                                                                                                                  0.00 
                                                                                                                                  \end{bmatrix}$ & $(1,1)$ \\                                                                                                                                                                                                                                                               
$\alpha^{fmin,j}_{lk}$       & $(1,1)$ & $\begin{bmatrix}
                                         0.99 \\
                                         0.01
                                        \end{bmatrix}$ & $\begin{bmatrix}
                                                         0.73 & 0.71\\
                                                         0.27 & 0.29
                                                      \end{bmatrix}$ & $(1,1)$  & $(1,1)$    & $\begin{bmatrix}
                                                                                               1.00\\
                                                                                               0.00
                                                                                                \end{bmatrix}$ & $\begin{bmatrix}
                                                                                                                 1.00\\
                                                                                                                 0.00
                                                                                                                 \end{bmatrix}$ & $\begin{bmatrix}
                                                                                                                                  1.00\\
                                                                                                                                  0.00 
                                                                                                                                  \end{bmatrix}$ & $(1,1)$ \\
 \hline                                                                                                                                 
$\beta^{2,j}_{k,l}$        & $\begin{bmatrix}
                             0.5\\
                             0.5 
                        \end{bmatrix}$  & $(1,1)$    & $\begin{bmatrix}
                                                       0.00 & 1.00\\
                                                       1.00 & 0.00
                                                       \end{bmatrix}$   & $\begin{bmatrix}
                                                                           0.00 \\
                                                                           1.00
                                                                           \end{bmatrix}$ & $\begin{bmatrix}
                                                                                       1.00 \\
                                                                                       0.00
                                                                                       \end{bmatrix}$ & $(1,1)$ & $(1,1)$ & $(1,1)$ & $\begin{bmatrix}
                                                                                                                                       0.5\\
                                                                                                                                        0.5
                                                                                                                                       \end{bmatrix}$\\                                                                                                                                                                                                                                                                            
$\alpha^{hyb,j}_{lk}$   & $(1,1)$ & $\begin{bmatrix}
                                         0.00 \\
                                         1.00
                                        \end{bmatrix}$ & $\begin{bmatrix}
                                                       0.00 & 1.00\\
                                                       1.00 & 0.00
                                                      \end{bmatrix}$ & $(1,1)$  & $(1,1)$    & $\begin{bmatrix}
                                                                                               0.00\\
                                                                                               1.00
                                                                                                \end{bmatrix}$ & $\begin{bmatrix}
                                                                                                                 1.00\\
                                                                                                                 0.00
                                                                                                                 \end{bmatrix}$ & $\begin{bmatrix}
                                                                                                                                  1.00\\
                                                                                                                                  0.00 
                                                                                                                                  \end{bmatrix}$ & $(1,1)$ \\
$\alpha^{ga,j}_{lk}$    & $(1,1)$ & $\begin{bmatrix}
                                          0.00 \\
                                          1.00
                                          \end{bmatrix}$ & $\begin{bmatrix}
                                                           0.00 & 1.00\\
                                                           1.00 & 0.0
                                                      \end{bmatrix}$ & $(1,1)$  & $(1,1)$    & $\begin{bmatrix}
                                                                                               0.00\\
                                                                                               1.00
                                                                                                \end{bmatrix}$ & $\begin{bmatrix}
                                                                                                                 1.00\\
                                                                                                                 0.00
                                                                                                                 \end{bmatrix}$ & $\begin{bmatrix}
                                                                                                                                  1.00\\
                                                                                                                                  0.00 
                                                                                                                                \end{bmatrix}$ & $(1,1)$ \\
 $\alpha^{fmin,j}_{lk}$  & $(1,1)$ & $\begin{bmatrix}
                                         0.00 \\
                                         1.00
                                        \end{bmatrix}$ & $\begin{bmatrix}
                                                         0.00 & 1.00\\
                                                         1.00 & 0.00
                                                      \end{bmatrix}$ & $(1,1)$  & $(1,1)$    & $\begin{bmatrix}
                                                                                               0.00\\
                                                                                               1.00
                                                                                                \end{bmatrix}$ & $\begin{bmatrix}
                                                                                                                 1.00\\
                                                                                                                 0.00
                                                                                                                 \end{bmatrix}$ & $\begin{bmatrix}
                                                                                                                                  1.00\\
                                                                                                                                  0.00 
                                                                                                                                  \end{bmatrix}$ & $(1,1)$ \\                                                                                                                                  
\hline                                                                                                                                  
\end{tabular}}
\label{TablaFollower}
\end{table}

\begin{table}
\tbl{Computational data of the different solvers for Case 2.}
{\begin{tabular}{|ccccc|}
\hline
Solver           & $J^\Delta_T (\alpha^{opt}, \beta^2)$ & iterations & evaluations & time (min) \\
\hline
\texttt{ga-hybrid}  & $1.1300\,  10^4$                  & 98    & 4955   & 57.93    \\
\texttt{ga}         & $1.1312\,  10^4$                  & 94    & 4750   & 47.84    \\
\texttt{fmincon}    & $1.1349\,  10^4$                  & 10    & 75    & 17.48   \\
\hline
\end{tabular}}
\label{Table_follower2}
\end{table}

\underline{Case 2:} For this second case we have chosen $\beta^{2,j}_{kl}$ such that we block the avenues $A_4$ and $A_8$ at intersections $j = 4, 5$ (that is, $\beta^{2,4}_{4\,12} = \beta^{2,5}_{8\,11} = 0$), and we let free pass of vehicles from $A_5$ to $A_7$ and also from $A_9$ to $A_6$ (that is, $\beta^{2,3}_{9\,6} = \beta^{2,3}_{5\,7} = 1$). Then, it is expected that drivers from $A_3$ and $A_{10}$ will turn on $A_5-A_7$ and $A_9-A_6$, respectively, avoiding the blocked avenues $A_8$ and $A_4$. Also for this case, drivers should take their outways at avenues $A_{12}$ and $A_{11}$, respectively. So, once fixed the vectors $\beta^{2,j}_{kl}$ in this manner (see full details in Table \ref{TablaFollower}), the follower problem was addressed again by using solvers \texttt{ga-hybrid}, \texttt{ga} and \texttt{fmincon-multi-start}. The corresponding preferences resulting from these solvers, denoted by $\alpha^{hyb,j}_{lk}, \alpha^{ga,j}_{lk}$ and $\alpha^{fmin,j}_{lk}$, can be also found in Table \ref{TablaFollower}.
 
For this case, the three solutions are practically equal, and the key of this fact relies in the achieved matrices of preferences at junction $j=3$. In these matrices, the values of preferences indicate that all drivers from $A_5$ prefer taking $A_7$ ($\alpha^3_{7\,5}=1$) and all drivers from $A_9$ prefer taking $A_6$ ($\alpha^3_{6\,9}=1$), avoiding in this way the blocked avenues $A_4$ and $A_8$. This could be checked in the level curves of pollution shown at Figure \ref{Caso2}. As in previous case, the predicted behavior of drivers is also fulfilled for this new case. 
  
Finally, the effectiveness and the computational cost (accuracy of solution, number of cost functional evaluations and computation time) of the three options to solve the follower problem were evaluated, but for the sake of conciseness only the output from Case 2 is shown. So, in Table \ref{Table_follower2} the discrete cost function $J^\Delta_T$ evaluated at the optimal solution, the number of functional evaluations and the solver execution time are shown. As can be seen there, the hybrid method and the genetic algorithm present a better effectiveness but at a much higher computational cost and execution time; at contrast, the \texttt{fmincon} solver has a slightly poorer effectiveness but with a significantly lower computational cost, being consistent with the smallest dimension of its initial vectors set.

\subsection{A highly restrictive Stackelberg solution}

\begin{figure}
\centering
\subfloat[Adjoint state.]{%
\resizebox*{5cm}{!}{\includegraphics[width = 10cm]{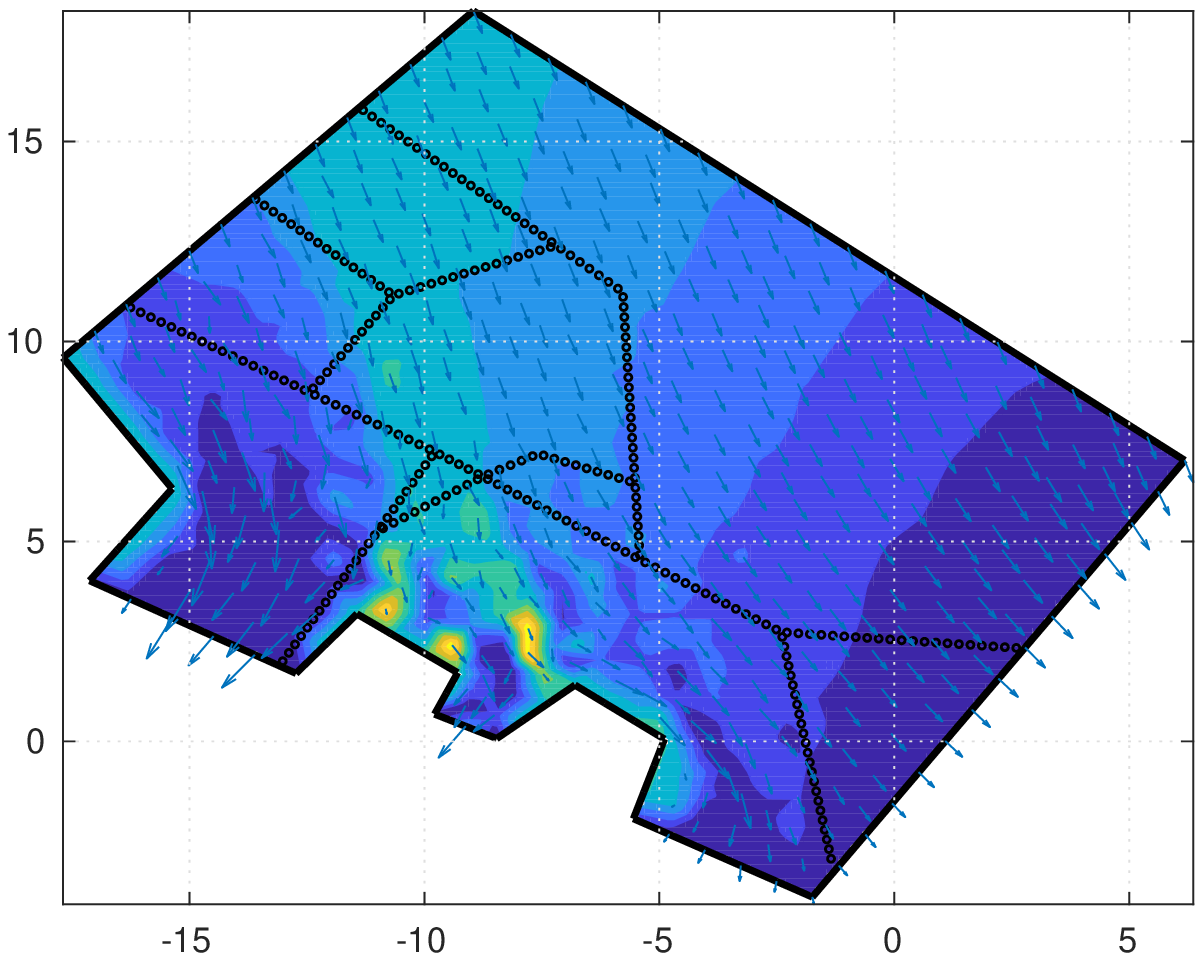}{\label{Isolineas_adjunto}}}}\hspace{5pt}
\subfloat[Number of iterations.]{%
\resizebox*{5cm}{!}{\includegraphics[width = 10cm]{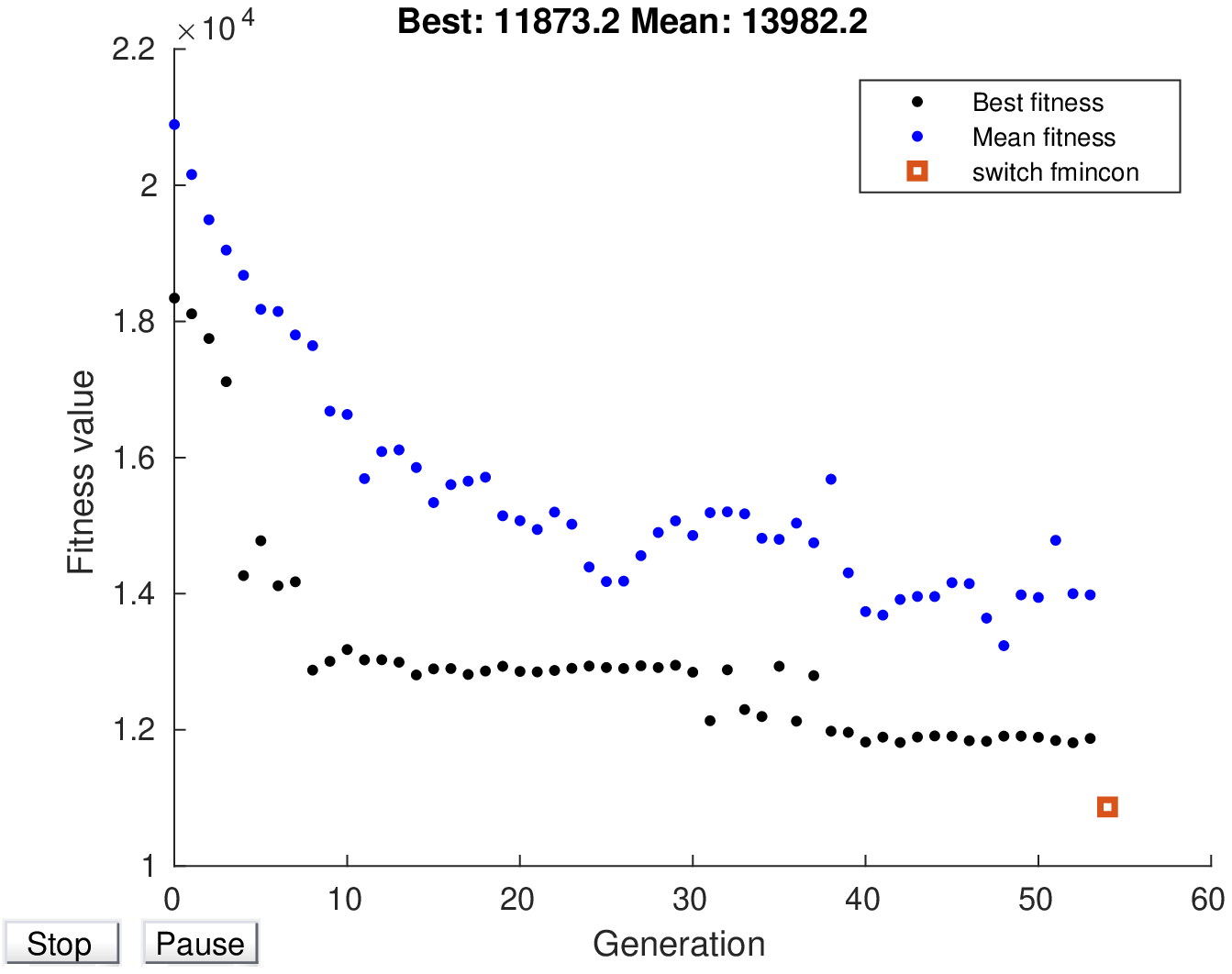}{\label{Iters_ga}}}}
\caption{(a) Level curves of the time-averaged adjoint discretized state employed for the optimization process. (b) Evolution of the number of iterations for the solver \texttt{ga-hybrid}, including the \texttt{fmincon} solution.}\label{Isolineas_NoOptimo1}
\end{figure}

As commented in above sections, as a previous step in the evaluation of $J^\Delta_P$, the adjoint model needs to be solved (only once). The isolines of the discretized adjoint state (averaged for a time interval of 24 hours) are shown in Fig. \ref{Isolineas_adjunto}. In this Figure, the adjoint state shows minimum values near the outflow boundary $S^-$ (except in a few small zones), and it increases as we get closer to the inflow boundary $S^+$. This tendency indicates that, in order to reduce the leader functional (which includes the product of the adjoint evaluated on roads and their corresponding emissions), the traffic must be concentrated in the roads located in the low zone of the domain. Therefore, it is expected that the leader's effort (using the restrictions at network intersections) should be directed to block the access of drivers to the top part of the domain.

The Stackelberg solution $(\alpha^*,\beta^*)$ is displayed at Table \ref{Tabla1}, and an analysis of the computational effort of the optimization process can be found at Fig. \ref{Iters_ga}. In this Stackelberg solution, restrictions $\beta^{*,3}_{5\,6} = \beta^{*,3}_{5\,7} = \beta^{*,4}_{4\,12} = 0$ indicate that the leader blocks the avenue $A_5$ at junction $j=3$, and avenue $A_4$ at junction $j=4$; on the contrary, $\beta^{*,5}_{8\,11}=1$ leaves free pass to cars from $A_8$ at junction $j=5$. With this strategy, the leader leaves to the follower with only an outway at $A_{11}-A_{14}$ avenues for those vehicles that enter by $A_{10}$, meanwhile vehicles that enter by $A_1$ and $A_2$ remains blocked in avenues $A_1,A_2,A_3,A_4$. This allows the leader to reach the objective of prevent drivers to take $A_6, A_{12}$ and $A_{15}$ which conduce to the upper part of the domain. With these blocked avenues the follower takes the (only) outway left for the leader ($\alpha^{*,6}_{8\,10} = 1$) and prefers the blocked $A_4$ instead of the blocked $A_5$ ($\alpha^{*,2}_{4\,3} = 1$). As a consequence of this preferences and restrictions, $(\alpha^*,\beta^*)$ are bottlenecks at intersections $j=1,2,4$, with large queue lengths at $A_1, A_2$ that maximize the length queue mean (Fig. \ref{fig7-c}), with low presence of vehicles at upper domain avenues (Fig. \ref{fig6-c}), but with a large mean density (Fig. \ref{fig7-a}), and low outflow at exit-points of the road network (Fig. \ref{fig7-d}). Therefore, a lower pollution levels are expected at most of the domain, and the higher levels are limited to specific zones close to the inferior boundary, as confirmed in Fig. \ref{fig8-c}. This situation is congruent with the prediction deduced above using the solution of the adjoint model.

It is important remarking here that, in contrast, the follower solution in above Case 2 (that will be taken here as a non-optimal case), do not present any blocked avenue and all densities are below $\rho_{max}$ in the whole network. Consequently, a better distribution of the density is reached (see Fig. \ref{fig6-a}) giving lower density, higher flow, less queues and finally more outflow with respect to the Stackelberg solution, but at expenses of a more polluted city (compare Figs. \ref{fig8-a} and \ref{fig8-c}).
 
All previous results suggest a worse situation for drivers in the Stackelberg solution, with large travel times to ensure low pollution levels. This is confirmed by the values displayed in Table \ref{Tabla2}, where the Stackelberg solution reduces the mean pollution by more than a $112\%$ but increases dramatically the travel time by four times.  

Is clear that this Stackelberg solution is easy of explain, agrees with all the methodology exposed and fulfils the predictions made. This strongly suggest that the correct solution was obtained, but in a realistic local government situation it is not possible to apply this (so restrictive) strategy to prevent the entry of vehicles in the city in order to reduce the urban pollution levels. So, in the following subsection we will present a Stackelberg solution with {\it relaxed} restrictions.

\subsection{A relaxed Stackelberg solution} 

In order to relax the traffic restrictions (avoiding the complete blockage of any avenue), the constraints (\ref{constraint2}) for the leader problem are changed to $0.2\leq \beta^j_{kl}\leq 0.8$ and $\sum_{k\in\mathcal{I}^{in}_j} \beta^j_{kl} = 1$. That is, the leader lets pass between $20\%$ and $80\%$ of vehicles from avenue $k\in\mathcal{I}^{in}_j$ to avenue $l\in\mathcal{I}^{out}_j$ at junction $j$, avoiding fully blocked avenues. 

\begin{figure}
\centering
\subfloat[Non-optimal case.]{%
\resizebox*{5cm}{!}{\includegraphics[width=4.5cm]{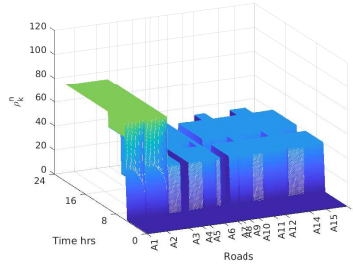}{\label{fig6-a}}}}\hspace{5pt}
\subfloat[Relaxed Stackelberg.]{%
\resizebox*{5cm}{!}{\includegraphics[width=4.5cm]{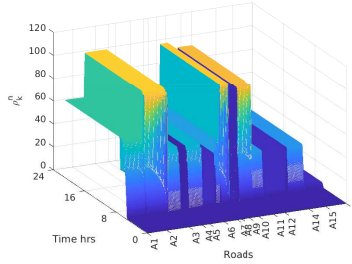}{\label{fig6-b}}}}
\subfloat[Restrictive Stackelberg.]{%
\resizebox*{5cm}{!}{\includegraphics[width=4.5cm]{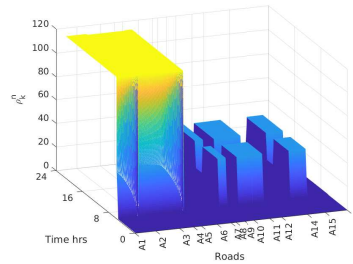}{\label{fig6-c}}}}
\caption{Time and spatial evolution of vehicle density at selected avenues for the three cases: non-optimal, relaxed Stackelberg and restrictive Stackelberg. The Stackelberg solutions generate congestion in some avenues ($\rho$ close to $\rho^{max} = 120 $).} 
\label{fig:fig6}
\end{figure}

\begin{figure}
\centering
\subfloat[Density]{%
\resizebox*{5cm}{!}{\includegraphics[width = 5 cm ]{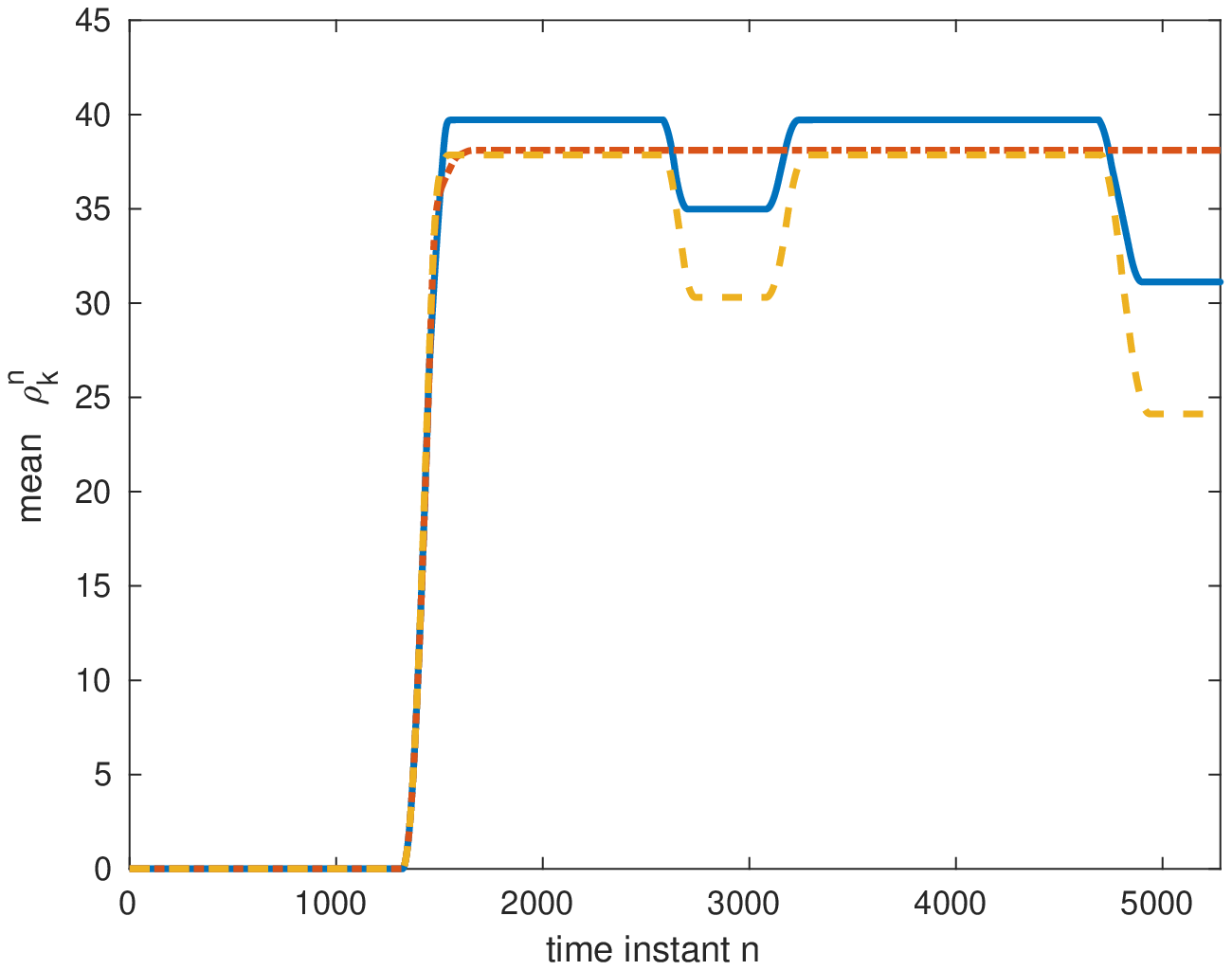}{\label{fig7-a}}}}\hspace{4pt}
\subfloat[Flow]{%
\resizebox*{5cm}{!}{\includegraphics[width = 5 cm ]{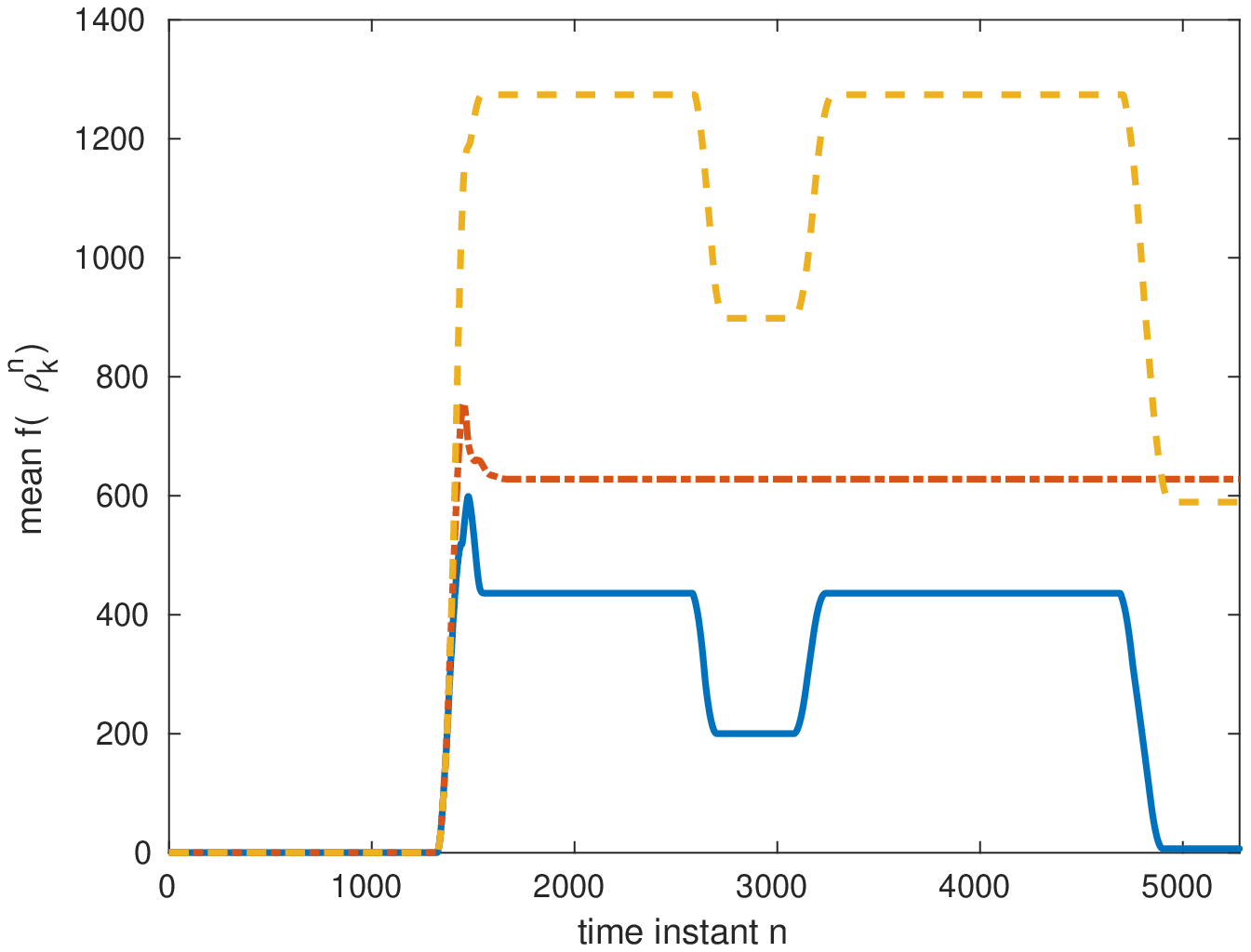}{\label{fig7-b}}}}\hspace{4pt}
\subfloat[Queue length]{%
\resizebox*{5cm}{!}{\includegraphics[width = 5 cm]{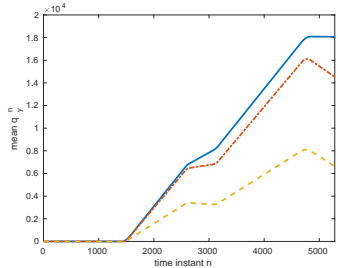}{\label{fig7-c}}}}\hspace{4pt}
\subfloat[Outflow]{%
\resizebox*{5cm}{!}{\includegraphics[width = 5 cm]{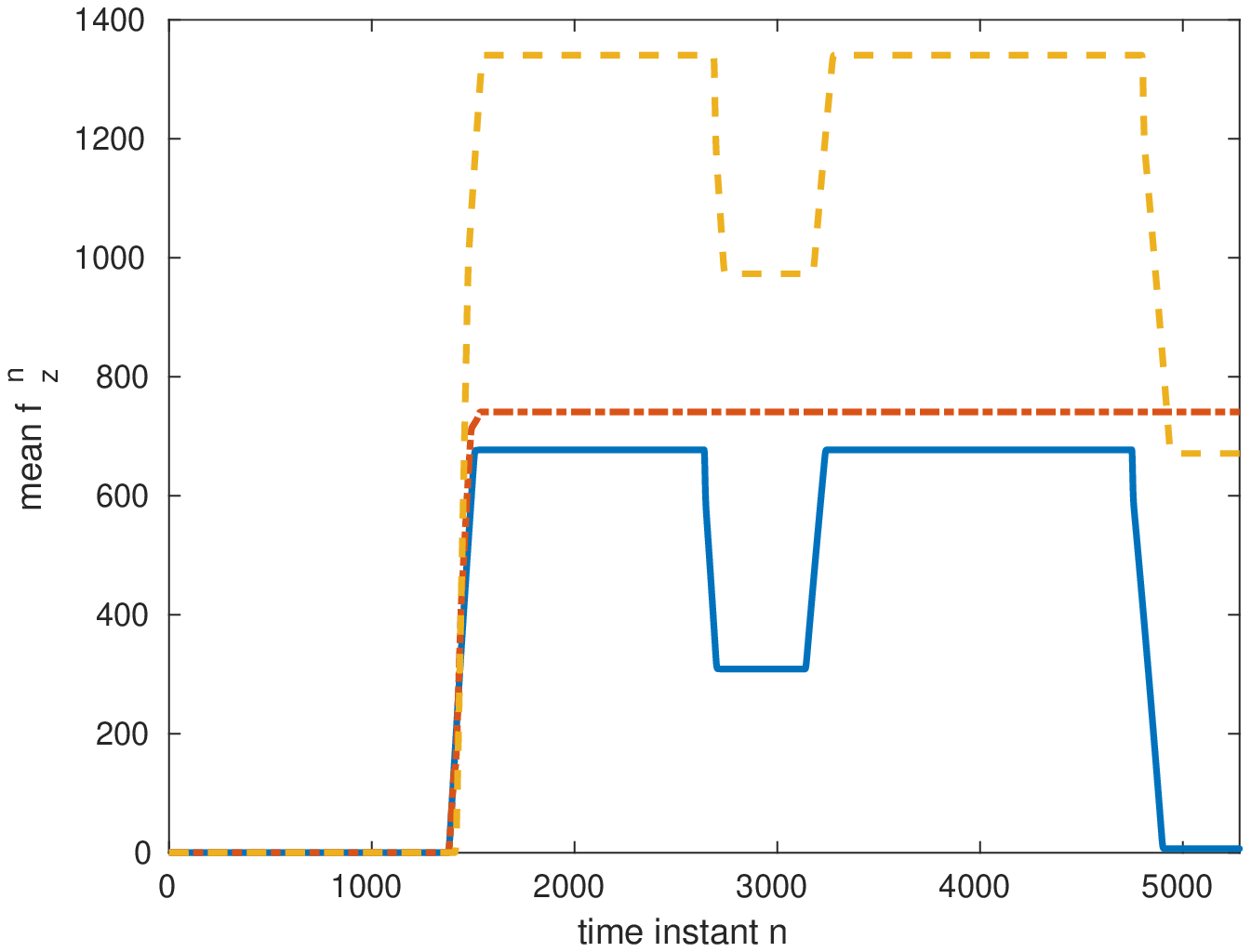}{\label{fig7-d}}}}
\caption{Mean values of traffic variables throughout the whole network corresponding to the highly restrictive Stackelberg (solid lines), the relaxed Stackelberg (dash-point lines), and the non-optimal (dashed lines).} \label{fig:fig7}
\end{figure}

\begin{figure}
\centering
\subfloat[Non-optimal case]{%
\resizebox*{5cm}{!}{\includegraphics[width=4.2cm]{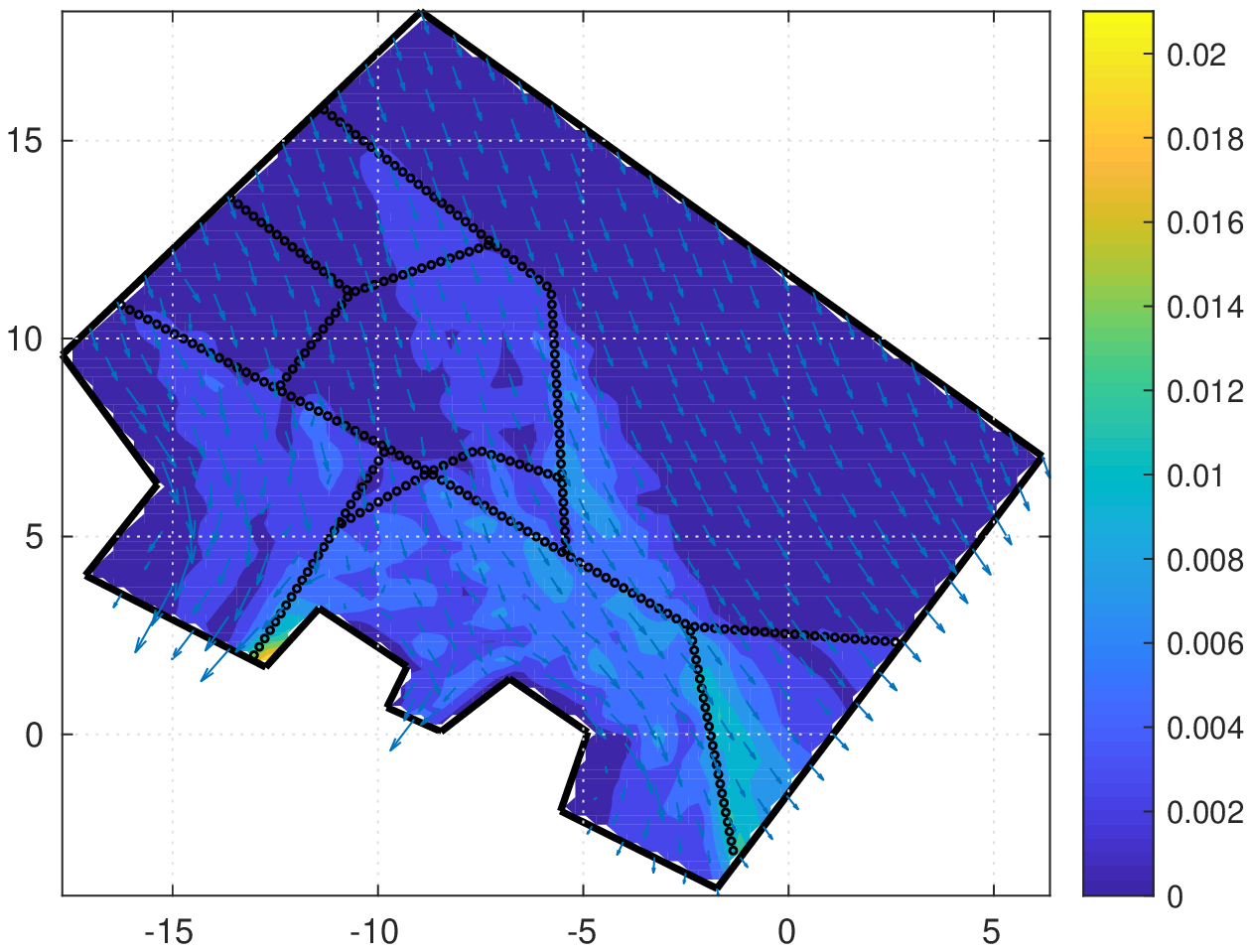}{\label{fig8-a}}}}\hspace{4pt}
\subfloat[Relaxed Stackelberg]{%
\resizebox*{5cm}{!}{\includegraphics[width=4.2cm]{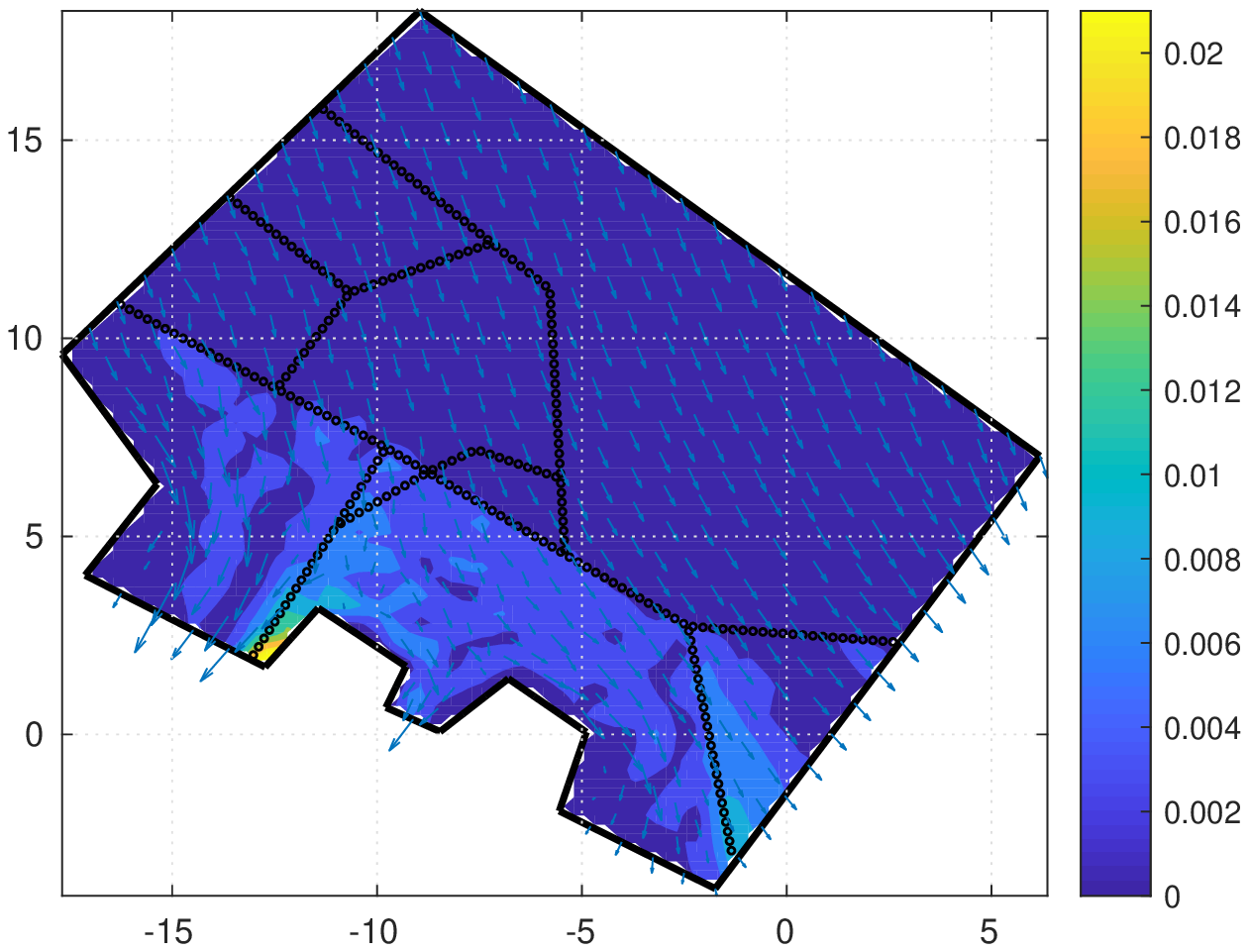}{\label{fig8-b}}}}\hspace{4pt}
\subfloat[Restrictive Stackelberg]{%
\resizebox*{5cm}{!}{\includegraphics[width=4.2cm]{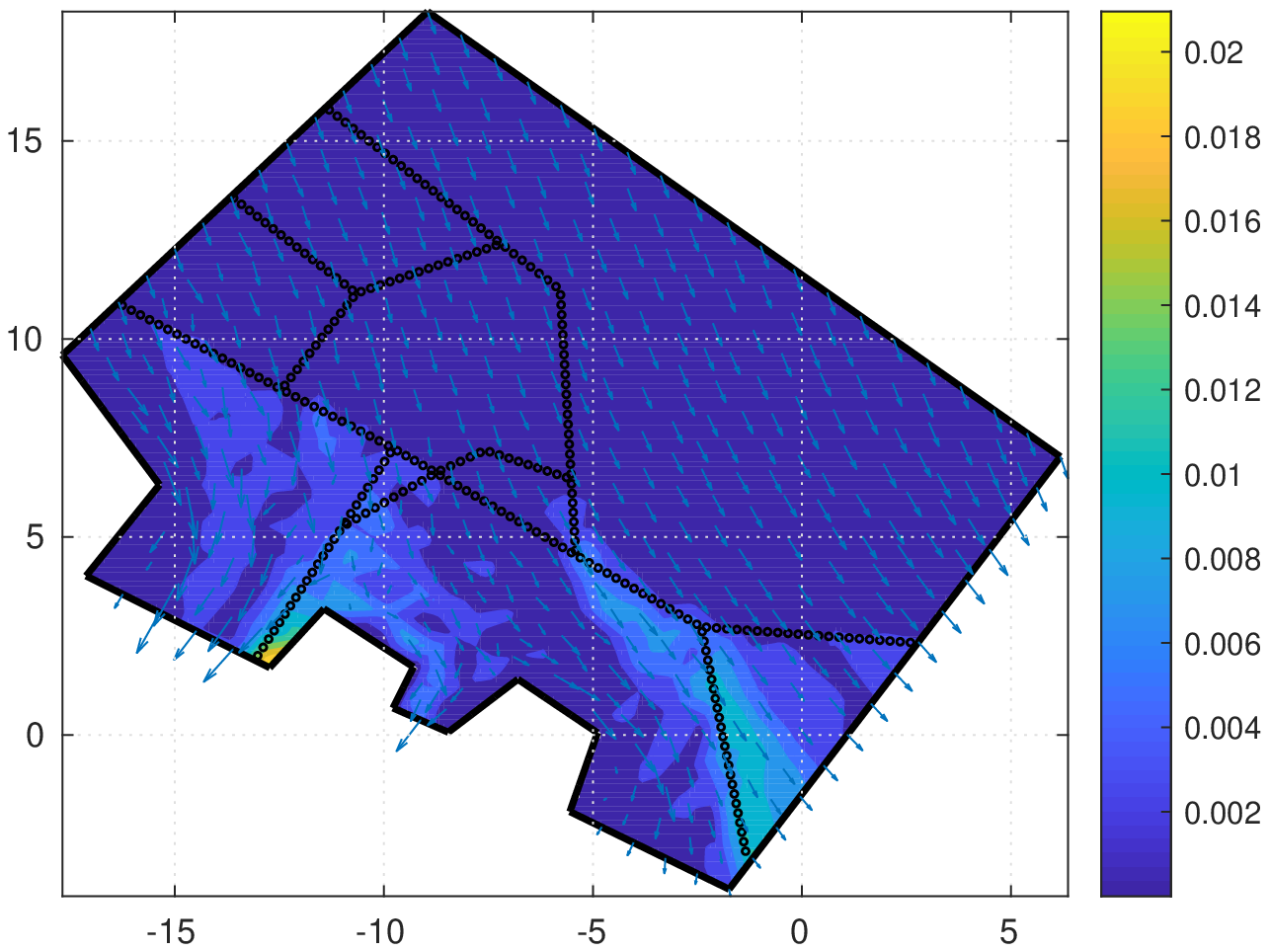}{\label{fig8-c}}}}\hspace{4pt}
\caption{Isolines for the time-averaged pollution concentration, obtained from the numerical resolution of model (\ref{MP}).}\label{fig:fig8}
\end{figure}

The relaxed Stackelberg solution $(\alpha^r,\beta^r)$ is qualitatively similar to the  restrictive case: It allows passing a minimum of vehicles to the upper part of the network, limits the traffic congestion to the inferior part of the city, and leaves to the follower the same outway by avenue $A_{14}$. This can be observed at Table \ref{Tabla1}, where the drivers' preferences $\alpha^{r,2}_{4,3} = 0, \alpha^{r,3}_{6\,5} = 0.08, \alpha^{r,3}_{6\,9} = 0.34, \alpha^{r,6}_{9\,10} = 0.06$ indicate that avenues $A_4,A_6,A_{12}$ are not chosen by drivers, allowing the leader to impose low restrictions $\beta^{r,1}_{9\,6}=\beta^{r,4}_{6\,12} = 0.20$. Meanwhile, the restrictions $\beta^{r,1}_{1\,3}= 0.79, \beta^{r,3}_{5\,7} = \beta^{r,5}_{11\,7} = 0.80$ give pass priority to the flows of vehicles at avenues $A_1,A_3,A_5,A_7,A_{11}$, limiting the traffic congestion to avenues $A_2,A_8,A_9$ in the inferior part of the network (see Fig. \ref{fig6-b} and Fig.\ref{fig8-b}).  

These restrictions and preferences generate the functional costs values given in Table \ref{Tabla2}, showing that this relaxed solution presents a decrease in its effectiveness (respect to the non-optimal Case 2) reducing the mean pollution in a $16\%$, but increasing travel time (although less than in the restrictive case). These results also suggest that the pollution levels increase as the restrictions are more relaxed, making harder for the leader to prevent traffic flows on the upper part of the network.

\begin{table}
\centering
\tbl{Values of Stackelberg solutions for the bi-level problem (\ref{P2f})-(\ref{P2l}): $(\alpha^*, \beta^*)$ stands for the highly restrictive case, and $(\alpha^r, \beta^r)$ for the relaxed one.}
{\begin{tabular}{|c|ccccccccc|}
      \hline
                           & j=1        & j=2            & j=3       & j=4       &  j=5      & j=6              & j=7         & j=8         & j=9 \\
         \hline
 $\mathcal{I}_j^{in}(k)$   & $\{1,2\}$  & $\{3\}$    & $\{5,9\}$ & $\{4,6\}$ & $\{7,8\}$ & $\{10\}$         & $\{12\}$    & $\{11\}$    & $\{16,17\}$\\
 $\mathcal{I}_j^{out}(l)$  & $\{3 \}$   & $\{4,5\}$  & $\{6,7 \}$& $\{12\}$  & $\{11\}$  & $\{8,9\}$        & $\{15,17\}$ & $\{14,16\}$ & $\{13\}$   \\
       \hline
 $\alpha^{*,j}_{lk}$        & $(1,1)$    & $\begin{bmatrix}
                                            1.00 \\
                                            0.00
                                      \end{bmatrix}$ & $\begin{bmatrix}
                                                       0.49 & 1.00\\
                                                       0.51 & 0.00
                                                      \end{bmatrix}$   & $(1,1)$ & $(1,1)$ & $\begin{bmatrix}
                                                                                            1.00\\
                                                                                            0.00
                                                                                            \end{bmatrix}$ & $\begin{bmatrix}
                                                                                                             1.00\\
                                                                                                             0.00
                                                                                                             \end{bmatrix}$ & $\begin{bmatrix}
                                                                                                                               1.00\\
                                                                                                                               0.00 
                                                                                                                              \end{bmatrix}$ & $(1,1)$ \\
$\beta^{*,j}_{kl}$     & $\begin{bmatrix}
                        0.00\\
                        1.00 
                        \end{bmatrix}$ & $(1,1)$ & $\begin{bmatrix}
                                                    0.00 & 0.00\\
                                                    1.00 & 1.00
                                                   \end{bmatrix}$   & $\begin{bmatrix}
                                                                     0.00 \\
                                                                     1.00
                                                                      \end{bmatrix}$ & $\begin{bmatrix}
                                                                                       0.00 \\
                                                                                       1.00
                                                                                       \end{bmatrix}$ & $(1,1)$ & $(1,1)$ & $(1,1)$ & $\begin{bmatrix}
                                                                                                                                        0.5\\
                                                                                                                                        0.5
                                                                                                                                       \end{bmatrix}$\\                                                                                                                                     
      \hline
$\alpha^{r,j}_{lk}$        & $(1,1)$    & $\begin{bmatrix}
                                            0.00 \\
                                            1.00
                                      \end{bmatrix}$ & $\begin{bmatrix}
                                                       0.08 & 0.34\\
                                                       0.92 & 0.66
                                                      \end{bmatrix}$   & $(1,1)$ & $(1,1)$ & $\begin{bmatrix}
                                                                                            0.94\\
                                                                                            0.06
                                                                                            \end{bmatrix}$ & $\begin{bmatrix}
                                                                                                             1.00\\
                                                                                                             0.00
                                                                                                             \end{bmatrix}$ & $\begin{bmatrix}
                                                                                                                               1.00\\
                                                                                                                               0.00 
                                                                                                                              \end{bmatrix}$ & $(1,1)$ \\
$\beta^{r,j}_{kl}$     & $\begin{bmatrix}
                        0.79\\
                        0.21
                        \end{bmatrix}$ & $(1,1)$ & $\begin{bmatrix}
                                                    0.79 & 0.79\\
                                                    0.21 & 0.21
                                                   \end{bmatrix}$   & $\begin{bmatrix}
                                                                     0.80 \\
                                                                     0.20
                                                                      \end{bmatrix}$ & $\begin{bmatrix}
                                                                                       0.80 \\
                                                                                       0.20
                                                                                       \end{bmatrix}$ & $(1,1)$ & $(1,1)$ & $(1,1)$ & $\begin{bmatrix}
                                                                                                                                        0.5\\
                                                                                                                                        0.5
                                                                                                                                       \end{bmatrix}$\\                                                                                                                                                                                                                                  
\hline                                                                                                                                             
\end{tabular}}
\label{Tabla1}
\end{table}

\begin{table}
\centering
\tbl{Numerical values of objective functionals $J^\Delta_P$ and $J^\Delta_T$ for different strategies.}
{\begin{tabular}{|ccc|}
     \hline
                         Case & $J^\Delta_P(\alpha,\beta)$ & $J^\Delta_T(\alpha,\beta)$\\
     \hline
     Restrictive Stackelberg  & $1.0866 \,10^4$  & $ 5.9695\,10^4$ \\ 
     Relaxed Stackelberg     & $1.9271 \,10^4$  & $ 4.8625\,10^4$ \\
     Non-optimal (Case 2)      & $2.2938 \,10^4$  & $ 1.1300\,10^4$ \\
     \hline     
\end{tabular}}
\label{Tabla2}
\end{table}        
         
\section{Conclusions}

In this work a bi-level optimal control problem was addressed in the sense of Stackelberg optimization. In the problem, the local government (leader) has as objective dropping down the urban pollution levels using traffic restrictions, meanwhile the drivers (follower) have the objective of minimizing their travel-time following their preferences. The optimal solution was obtained using a combination of genetic and interior-point algorithms applied to previous numerical simulations of the traffic on an urban road network, of its pollution emissions and of the pollutant transport over the whole urban domain.    

With the aim of saving computational efforts, adjoint techniques were used. However, the impossibility of deriving an exact gradient of the objectives with respect to the restrictions and preferences implied a large computational cost for getting a Stackelberg solution. 

The numerical experiences showed that the effectiveness of the Stackelberg solution is higher when the restrictions were such that complete blockage of traffic at road intersection is allowed, minimizing the pollution levels and increasing the travel time. This effectiveness presents a significant drop down when the restrictions are relaxed. This fact agrees, in a roughly way, with the above referenced empirical studies where data shown a progressive increase of pollution levels from non-traffic zones to traffic ones.  

Finally, future research work could be related with modifications in the objectives for both, the leader and the follower. So, the queue length could be considered as an additional leader objective, removing it from the follower cost. Also, future work could consider more sophisticated improvements in the traffic model, making it more complex and realistic for urban domains. The macroscopic models with dynamic velocity and a LWR model with diffusion and different forcing (traffic-lights, multiple lanes, in-out ramps\dots) are options available in the specialized literature \cite{treiber}.

\section*{Acknowledgements}

This work was supported by Ministerio de Econom\'{\i}a y Competitividad (Spain)/FEDER under Grant MTM2015-65570-P; Xunta de Galicia under Grant ED431C 2018/50.
The first author also thanks the support from Sistema Nacional de Investigadores (Mexico) under Grant SNI-52768; Programa para el Desarrollo Profesional Docente (Mexico) under Grant PRODEP/103.5/16/8066.

\section*{Disclosure statement}
The authors declare no conflict of interest with the paper

\end{document}